%% file: so35.tex
\documentclass{amsart}

\newtheorem{theorem}{Theorem}[section]
\newtheorem{lemma}[theorem]{Lemma}

\newtheorem{proposition}[theorem]{Proposition}
\newtheorem{corollary}[theorem]{Corollary}

\theoremstyle{definition}
\newtheorem{definition}[theorem]{Definition}

\theoremstyle{remark}
\newtheorem{remark}[theorem]{Remark}

\usepackage{graphicx}
\usepackage{color}
\usepackage{amsmath}
\usepackage{amsfonts}
\usepackage{amssymb}
\newcommand{\be}{\begin{equation}}

\newcommand{\ee}{\end{equation}}

\newcommand{\al}{\alpha}
\newcommand{\alt}{\tilde{\alpha}}

\newcommand{\Om}{\Omega}

\newcommand{\lc}{{\stackrel{\scriptscriptstyle{LC}}{\Gamma}}}

\newcommand{\dz}{\wedge}

\newcommand{\ba}{\begin{array}}

\newcommand{\ea}{\end{array}}

\newcommand{\beq}{\begin{eqnarray}}
\newcommand{\riclc}{Ric^{\scriptscriptstyle{LC}}}
\newcommand{\ricga}{Ric^{\scriptscriptstyle{\Gamma}}}
\newcommand{\eeq}{\end{eqnarray}}

\newtheorem{lm}{lemma}

\newtheorem{thee}{theorem}

\newtheorem{proo}{proposition}

\newtheorem{co}{corollary}

\newtheorem{rem}{remark}

\newtheorem{deff}{definition}

\newcommand{\bd}{\begin{deff}}

\newcommand{\ed}{\end{deff}}

\newcommand{\bl}{\begin{lm}}

\newcommand{\el}{\end{lm}}

\newcommand{\bp}{\begin{proo}}

\newcommand{\ep}{\end{proo}}

\newcommand{\bt}{\begin{thee}}

\newcommand{\et}{\end{thee}}

\newcommand{\bc}{\begin{co}}

\newcommand{\ec}{\end{co}}

\newcommand{\brm}{\begin{rem}}

\newcommand{\erm}{\end{rem}}

\newcommand{\der}{{\rm d}}

\newcommand{\frg}{\frak{g}}

\newcommand{\ett}{\tilde{\eta}}
\newcommand{\sgn}{\mathrm{sgn}}

\usepackage{t1enc}
\def\frak{\mathfrak}

\newcommand{\scm}{{\scriptscriptstyle{-}}}

\newcommand{\newc}{\newcommand}

\newcommand{\Ad}{\operatorname{Ad}}

\newcommand{\id}{\operatorname{id}}

\let\ccdot\cdot
\def\cdot{\hbox to 2.5pt{\hss$\ccdot$\hss}}

\newc{\aR}{\mbox{\boldmath{$ R$}}}
\newc{\aS}{\mbox{\boldmath{$ S$}}}
\newc{\aT}{\mbox{\boldmath{$ T$}}}
\newc{\aW}{\mbox{\boldmath{$ W$}}}

\newc{\aK}{\mbox{\boldmath{$ K$}}}
\newc{\aL}{\mbox{\boldmath{$ L$}}}

\newcommand{\twist}{\mathbb{T}}
\newcommand{\bbP}{\mathbb{P}}
\newcommand{\bbC}{\mathbb{C}}
\newcommand{\calV}{\mathcal{V}}
\newcommand{\calH}{\mathcal{H}}
\newcommand{\sph}{\mathbb{S}}

\usepackage{amssymb}
\usepackage{amscd}

\newcommand{\End}{\operatorname{End}}

\def\bbZ{{\mathbb{Z}}}

\newcommand{\bas}{{\bf e}}

\newc{\obstrn}[2]{B^{#1}_{#2}}

\newcommand{\rpl}                         
{\mbox{$
\begin{picture}(12.7,8)(-.5,-1)
\put(0,0.2){$+$}
\put(4.2,2.8){\oval(8,8)[r]}
\end{picture}$}}

\newcommand{\lpl}                         
{\mbox{$
\begin{picture}(12.7,8)(-.5,-1)
\put(2,0.2){$+$}
\put(6.2,2.8){\oval(8,8)[l]}
\end{picture}$}}

\usepackage{ifthen}
\newcommand{\bbR}{\mathbb{R}}
\newcommand{\sog}{\mathbf{SO}}
\newcommand{\spin}{\frak{spin}}
\newcommand{\slg}{\mathbf{SL}}
\newcommand{\og}{\mathbf{O}}
\newcommand{\soa}{\frak{so}}

\newcommand{\sua}{\frak{su}}
\newcommand{\ua}{\frak{u}}
\newcommand{\dr}{\mathrm{d}}
\newcommand{\sug}{\mathbf{SU}}
\newcommand{\ug}{\mathbf{U}}
\newcommand{\gat}{\tilde{\gamma}}
\newcommand{\Gat}{\tilde{\Gamma}}
\newcommand{\thet}{\tilde{\theta}}
\newcommand{\Thet}{\tilde{T}}

\newcommand{\st}{\sqrt{3}}
\newcommand{\kat}{\tilde{\kappa}}

\newc{\tensor}[1]{#1}
\newc{\Mvariable}[1]{\mbox{#1}}
\newc{\down}[1]{{}_{#1}}
\newc{\up}[1]{{}^{#1}}

\newc{\JulyStrut}{\rule{0mm}{6mm}}
\newc{\midtenPan}{\mbox{\sf S}}
\newc{\midten}{\mbox{\sf T}}
\newc{\midtenEi}{\mbox{\sf U}}
\newc{\ATen}{\mbox{\sf E}}
\newc{\BTen}{\mbox{\sf F}}
\newc{\CTen}{\mbox{\sf G}}

\def\sideremark#1{\ifvmode\leavevmode\fi\vadjust{\vbox to0pt{\vss
 \hbox to 0pt{\hskip\hsize\hskip1em
 \vbox{\hsize3cm\tiny\raggedright\pretolerance10000
 \noindent #1\hfill}\hss}\vbox to8pt{\vfil}\vss}}}%
                        
\newcommand{\bgw}{{\textstyle \bigwedge}}
\newcommand{\bgs}{{\textstyle \bigodot}}
\newcommand{\bgt}{{\textstyle \bigotimes}}
\newcommand{\Span}{\mathrm{Span}}

\numberwithin{equation}{section}

\newcounter{romenumi}
\newcommand{\labelromenumi}{(\roman{romenumi})}

\newenvironment{romlist}{\begin{list}{\labelromenumi}{\usecounter{romenumi}}}{\end{list}}

\newcommand{\ten}{\Upsilon}
\newcommand{\ft}{\curlyvee}

\begin{document}
\title{Irreducible SO(3) geometry in dimension five} 
\vskip 1.truecm \author{Marcin Bobie\'nski} \address{Instytut Matematyki,
Universytet Warszawski, ul. Banacha 2, Warszawa, Poland}
\email{mbobi@mimuw.edu.pl} \thanks{This research was supported by
the KBN grants 1 P03A 01529 and 1 P03B 07529}

\author{Pawe\l~ Nurowski} \address{Instytut Fizyki Teoretycznej,
Uniwersytet Warszawski, ul. Hoza 69, Warszawa, Poland}
\email{nurowski@fuw.edu.pl} \thanks{During the preparation of this
article the authors were members of the VW Junior Research Group ``Special
Geometries in Mathematical Physiscs'' at Humboldt University in
Berlin.}

\date{\today}

\begin{abstract}  
We consider the nonstandard inclusion of $\sog(3)$ in $\sog(5)$ associated
with a 5-dimensional irreducible representation. The tensor $\ten$ 
representing this reduction is found to be given by a ternary
symmetric form with special properties. A 5-dimensional manifold
$(M,g,\ten)$ with Riemannian metric $g$ and ternary form generated 
by such a tensor has a corresponding $\sog(3)$ structure, whose 
Gray-Hervella type classification is established using
$\soa(3)$-valued connections with torsion.

Structures with antisymmetric torsions, we call them the nearly
integrable $\sog(3)$ structures, are studied in detail. 
In particular, it is shown that the integrable models 
(those with vanishing torsion) are isometric to the symmetric spaces
$M_+= \sug(3)/\sog(3)$, $M_-=\slg(3,R)/\sog(3)$, $M_0=\bbR^5$.
We also find all nearly integrable $\sog(3)$ structures with
transitive symmetry groups of dimension $d>5$ and some examples for
which $d=5$.

Given an $\sog(3)$ structure $(M,g,\ten)$, we define its "twistor
space" $\twist$ to be the $\sph^2$-bundle of those unit 2-forms on $M$
which span $\bbR^3=\soa(3)$. The 7-dimensional twistor manifold
$\twist$ is then naturally equipped with several $CR$ and $G_2$
structures. The ensuing integrability conditions are discussed 
and interpreted in terms of the Gray-Hervella type classification. 

\vskip5pt\centerline{\small\textbf{MSC classification}: 53A40, 53B15,
  53C10}\vskip15pt
\end{abstract}
\maketitle
\tableofcontents
\section{Introduction}
In Cartan's list of the irreducible symmetric spaces of Type I 
the first entry is occupied by the family of symmetric spaces
$\sug(n)/\sog(n)$. If $n=2$ the
corresponding manifold is a 2-dimensional sphere $\sph^2$, but $n=3$
already corresponds to a nontrivial manifold
$M_+=\sug(3)/\sog(3)$. This is the so called Wu space \cite{Sm,Wuuuu}
which has a number of interesting
properties. Among them there is a fact that $M_+$ constitutes the
lowest dimensional example of a simply
connected manifold {\it not} admitting a $\mathbf{Spin}^c$ structure
\cite{ML}. From the point of view of the present paper another
property of this space is crucial: the isotropy representation 
of $M_+=\sug(3)/\sog(3)$ coincides with the irreducible 5-dimensional 
representation of $\sog(3)$. Thus, this space provides a symmetric 
model of a 5-dimensional manifold equipped with the irreducible
$\sog(3)$ structure. Inspecting the entire 
Cartan list of the irreducible symmetric spaces one finds (in Type
III, again at  the first entry!) another 5-dimensional space  
$M_-=\slg(3,\bbR)/\sog(3)$ equipped with the natural irreducible $\sog(3)$
structure.

The aim of this paper is to study 5-dimensional geometries modelled on 
the spaces $M_+$ and $M_-$. By this we mean studies of 5-dimensional 
manifolds with the reduction of the structure group of the 
$\sog(5)$-frame bundle to the irreducible $\sog(3)$. This
places the paper in the domain of {\it special geometries}, 
i.e. Riemannian geometries equipped with additional geometric
structures. In Ref. \cite{Fried} Th. Friedrich provides a general framework for
analysing such geometries. 
He also proposes the investigation of geometries modelled
on $M_+$ there. 

The framework for analysis of special geometries consists of
several steps. First, one distinguishes a
geometric object, preferably of tensorial type, that reduces the structure 
group and defining the special geometry. Then, one introduces a
metric connection which preserves this
object. As the last step one determines the restrictions on the
special geometry for this connection to be unique. This unique
connection, its torsion and curvature are then the main tools to study
the properties of the considered special geometry. 

It is instructive to illustrate this
procedure on the well known example of a {\it nearly} K\"ahler
geometry. Our choice of nearly K\"ahler geometry for this illustration is
motivated by the fact that its behaviour is remarkably close \cite{PA}
to all the phenomena we want to discuss in the context of the irreducible 
$\sog(3)$ geometries in dimension five. 

A Riemannian geometry $(M,g)$ on a $2n$-dimensional manifold $M$ can be made
more special by an introduction of a metric compatible almost complex structure. This is a tensor field $J:TM\to TM$ which
satisfies $J^2=-\id$ and 
$g(JX,JY)=g(X,Y)$. The tensor $J$ reduces the structure group from
$\sog(2n)$ to $\ug(n)$ and induces the distinguished inclusion
of the Lie algebra $\ua(n)$ in $\soa(2n)$. This inclusion
defines a class of a metric compatible connections $\Gamma$ which preserve
$J$. Here and in the following we will represent connections by means
of Lie-algebra-valued 1-forms on manifolds so, in the considered
case, $\Gamma\in\ua(n)\otimes\Omega^1(M)$, where
$\ua(n)\subset\soa(2n)$. The connections $\Gamma$ are highly
not unique. However, since all of them may be considered as elements 
of $\soa(2n)\otimes\Omega^1(M)$, i.e. as elements of the space in which the
Levi-Civita connection $\lc$ resides, one can try to make $\Gamma$
unique by the requirement that in the decomposition 
\be
\lc=\Gamma+\tfrac12 T\label{intr1}
\ee
the $T$-part has some special properties. In the considered case the 
uniqueness of $\Gamma$ is achieved by the demand that in the above
decomposition 
\be
\quad T\in\Omega^3(M).\label{intr2}
\ee
The 3-form $T$ is then interpreted as a skew-symmetric torsion of
the connection $\Gamma$. It follows that the decomposition 
(\ref{intr1})-(\ref{intr2}) is possible only for a subclass of metric
compatible almost complex structures. They may be characterised by the
condition 
$$(\stackrel{\scriptscriptstyle{LC}}{\nabla}_vJ)(v)=0\quad\quad\forall v\in TM.$$
The metric compatible almost complex structures satisfying this
condition are called {\it nearly} K\"ahler. Their geometric properties
are described in terms of the properties of the unique $\ua(n)$-valued
connection $\Gamma$ defined by (\ref{intr1})-(\ref{intr2}). In
particular, the torsion-free case, $T\equiv 0$, corresponds to
K\"ahler geometries. Another types of the nearly K\"ahler structures
may be distinguished by specifying that the curvature of $\Gamma$
belongs to a particular $\ug(n)$-irreducible
component of the tensor representation $\ua(n)\otimes\Omega^2(M)$.

Our treatment of the irreducible $\sog(3)$ geometries in dimension
five imitates the above approach to the nearly-K\"ahler geometries. We
first introduce an object, the (3,0)-rank tensor $\ten$, which reduces
the $\sog(5)$ structure to the irreducible $\sog(3)$. Although this
tensor has a different rank then $J$ its geometric characterisation,
which is a certain algebraic quadratic identity on $\ten$,
resembles very much the quadratic condition $J^2=-\id$. Using $\ten$
we distinguish an inclusion of $\soa(3)$ in $\soa(5)$. This maximal
inclusion is used on a Riemannian manifold endowed with
$\ten$ to distinguish a class of $\soa(3)$-valued metric
connections $\Gamma$. These are such that, in the decomposition
(\ref{intr1}), they have the skew-symmetric $T$-part. 
It follows that such connections, if exist, are unique. Their
existence is only possible for a particular class of tensors $\ten$
characterised by the condition
$$(\stackrel{\scriptscriptstyle{LC}}{\nabla}_v\ten)(v,v,v)=0\quad\quad\forall
v\in TM.$$

The organisation of the paper is reflected in the table of
contents. The notation is standard. However, 
depending on the context and esthetics of the presentation, we use
both the Schouten notation with the indices of tensors as well as the
geometric, index-free notation. Since all the time we are in the
Riemannian category, we do not distinguish between covariant and
contravariant tensors. This convention, when used in the formulae
written in the Schouten notation, enables as to identify tensors with
upper and lower indices. We will write them in the both positions
depending on convenience. In the entire text the Einstein summation
convention is assumed.

\section{Tensor $\ten$ reducing ${\bf O}(5)$ to the irreducible ${\bf
    SO}(3)$}
The two obvious examples $M_+={\bf SU}(3)/{\bf SO}(3)$ and $M_-={\bf
  SL}(3,\mathbb{R})/{\bf SO}(3)$ of the irreducible ${\bf SO}(3)$
structures should be supplemented by still another one, which in a
certain sense, is the simplest. One achieves this example 
by identifying vectors $A$ in 
$\mathbb{R}^5$ with $3\times 3$ symmetric traceless real
matrices $\sigma(A)$, 
\be
\mathbb{M}^5=\{~\sigma(A)\in M_{3\times 3}(\mathbb{R}):~~
\sigma(A)^{\rm T}=\sigma(A),~~{\rm tr}(\sigma(A))=0~\},\label{r5}
\ee
and defining the unique irreducible 5-dimensional representation
$\rho$ of
${\bf SO}(3)$ in $\bbR^5$ by 
\be
\rho(h)A~=~h~\sigma(A)~h^{\rm T},\quad\quad\quad\forall~h\in{\bf SO}(3),~~~~~~~A\in \bbR^5. \label{so35}
\ee
Then  
$M_0=({\bf SO}(3)\times_\rho\mathbb{R}^5)/{\bf SO}(3)$ also has an
irreducible $\sog(3)$ structure. 

From now on we identify $\mathbb{R}^5$ with matrices $\mathbb{M}^5$ as in
(\ref{r5}). Given an element $A\in\mathbb{R}^5$ we consider its characteristic
polynomial 
$$
P_A(\lambda)={\rm det}(\sigma(A)-\lambda I)=-\lambda^3+g(A,A)\lambda+\frac{2\sqrt{3}}{9}\ten(A,A,A).
$$
This polynomial is invariant under the ${\bf SO}(3)$-action given by
the representation $\rho$ of (\ref{so35}), 
$$
P_{\rho(h)A}(\lambda)=P_A(\lambda).
$$
Thus, all the coefficients
of $P_{A}(\lambda)$, which are multilinear in $A$, are 
${\bf SO}(3)$-invariant. It is convenient to choose a basis $\bas_i$ in
$\bbR^5$ in such a way that the identification $\sigma$ is given by 
\be
\bbR^5\ni A=a_i\bas_i\longmapsto \sigma(A)=\begin{pmatrix} 
\frac{a^1}{\sqrt{3}}-a^4&a^2&a^3\\
a^2&\frac{a^1}{\sqrt{3}}+a^4&a^5\\
a^3&a^5&-2\frac{a^1}{\sqrt{3}}
\end{pmatrix}\in\mathbb{M}^5.\label{bass}
\ee
After this convenient choice, the bilinear form $g$ simply becomes
\be
g(A,A)=a_1^2+a_2^2+a_3^2+a_4^2+a_5^2\label{mee},
\ee
and the ternary one $\ten$ is given by
\be
\ten(A,A,A)=
\frac{1}{2}a_1\big(~6a_2^2+6a_4^2-2a_1^2-3a_3^2-3a_5^2~\big)+\frac{3\sqrt{3}}{2}a_4(a_5^2-a_3^2)+3\sqrt{3}a_2a_3a_5.\label{tensor}
\ee
Both $g$ and $\ten$ are obviously ${\bf SO}(3)$-invariant. Since $g$ 
is the usual Riemannian metric on $\bbR^5$ the action
$\rho$ of (\ref{so35}) gives a nonstandard {\it irreducible} inclusion  
\be 
\iota: \sog(3)\hookrightarrow\og(5). \label{iota}
\ee
\begin{remark}
Although it is obvious we remark that
$$\ten(A,A,A)=\frac{3\sqrt{3}}{2}\det(\sigma(A)).$$
\end{remark}
In the following we consider a tensor $\ten_{ijk}\in\bigodot^3\bbR^5$ such that
$$\ten(A,A,A)=\ten_{ijk}a^ia^ja^k.$$
A simple algebra leads to the following proposition.
\begin{proposition}
\label{pr:ident}
The tensor $\ten_{ijk}$ has the following properties
\begin{itemize}
\item[i)] it is totally symmetric, $\ten_{ijk}=\ten_{(ijk)}$,
\vskip1mm
\item[ii)]it is trace-free, $\ten_{ijj}=0$,
\vskip1mm
\item[iii)] it satisfies the following identity
$$\ten_{jki}\ten_{lni}+\ten_{lji}\ten_{kni}+\ten_{kli}\ten_{jni}=g_{jk}g_{ln}+g_{lj}g_{kn}+g_{kl}g_{jn},$$ 
where
  $g(A,A)=g_{ij}a^i a^j$.
\end{itemize}
\end{proposition}
\begin{remark}
It is worth noting that property iii) after contraction with $g_{kn}$ and
$\ten_{mkn}$, respectively, implies
\begin{eqnarray*}
4 \ten_{ijk}\ten_{mjk}&=&14g_{im}, \\
4 \ten_{ilm}\ten_{jln}\ten_{kmn}&=&-3\ten_{ijk}.
\end{eqnarray*}
\end{remark}
Group ${\bf O}(5)$ naturally acts on $\bigodot^3 \mathbb{R}^5$ by
$$
\ten_{ijk}\mapsto H^{l}_{~i}H^{m}_{~~j}H^{n}_{~k}\ten_{lmn}, \quad\quad H\in\og(5).
$$
Our aim now is to find the stabiliser $G_\ten$ of tensor $\ten_{ijk}$ under
this action. We know that ${\bf SO}(3)\subset G_\ten$. In
the following we show that it is actually equal to ${\bf SO}(3)$. To
see this we take a 1-parameter subgroup $H(s)={\rm e}^{sX}$ of ${\bf
  SO}(5)$ generated by an element $X$ of the Lie algebra
$\frak{so}(5)$ in the standard 5-dimensional representation of skew
symmetric matrices. Taking $\frac{\rm d}{\rm ds}_{|s=0}$ of the
stabilising equation
$
\ten_{ijk}=H(s)^{l}_{~i}H(s)^{m}_{~~j}H(s)^{n}_{~k}\ten_{lmn}
$
we get the following linear equation 
\be
\ten_{ljk}X^{l}_{~~i}+\ten_{ilk}X^{l}_{~~j}+\ten_{ijl}X^{l}_{~~k}=0\label{defso3}
\ee
for the elements of the Lie algebra of the stabiliser. Its general
solution is 
$$
X=(X^i_{~j})=x^1E_1+x^2E_2+x^3E_3=x^I E_I,
$$ 
where $(x^I)$, $I=1,2,3,$ are real parameters and the matrices 
\be
E_1=\left(\begin{smallmatrix} 
0&0&0&0&\scriptscriptstyle{\sqrt{3}}\\
0&0&1&0&0\\
0&\scm1&0&0&0\\
0&0&0&0&1\\
\scriptscriptstyle{-\sqrt{3}}&0&0&\scm1&0
\end{smallmatrix}\right),\quad\quad\quad E_2=\left(\begin{smallmatrix} 
0&0&\scriptscriptstyle{\sqrt{3}}&0&0\\
0&0&0&0&1\\
\scriptscriptstyle{-\sqrt{3}}&0&0&1&0\\
0&0&\scm1&0&0\\
0&\scm1&0&0&0
\end{smallmatrix}\right),\quad\quad
E_3=\left(\begin{smallmatrix} 
0&0&0&0&0\\
0&0&0&2&0\\
0&0&0&0&1\\
0&\scm2&0&0&0\\
0&0&\scm1&0&0
\end{smallmatrix}\right),\label{so3bas}
\ee
satisfy the $\frak{so}(3)$ commutation
relations 
$$
[E_1,E_2]=E_3,\quad\quad[E_3,E_1]=E_2,\quad\quad
[E_2,E_3]=E_1,
$$
or $[E_J,E_K]=\epsilon^I_{JK}E_I$, for short.
Thus, the intersection of the stabiliser with the $\sog(5)$ component of $\og(5)$ is equal to the irreducible $\sog(3)$. 
Actually the stabiliser does not intersect with the complement of $\sog(5)$ in $\og(5)$, as it is explained in the following lemma.
\begin{lemma}
\label{lem:ginso}
The stabiliser of $\ten_{ijk}$ is contained in $\sog(5)$ component of $\og(5)$.
\end{lemma}
\begin{proof}
Since the complement of $\sog(5)$ in $\og(5)$ consists of elements of
the form $-g$ such that $g\in\sog(5)$ it is enough to prove that $-g$
with $g\in\sog(5)$ can not be in $G_\ten$.  
Assuming the opposite i.e. that $g\in\sog(5)$ and $-g\in G_\ten$  we get
the contradiction by the following steps.
The adjoint map $\Ad_g$ 
preserves $\soa(3)$. Thus it provides an orthogonal (with respect to the Killing form) transformation of $\soa(3)$
\[
\Ad_g|_{\soa(3)} \in \sog(\soa(3)),\qquad \soa(3) = \Span(E_1,E_2,E_3).
\]
On the other hand, any orthogonal transformation of our $\soa(3)$ has
the form $\Ad_h$ for an element $h\in\iota(\sog(3))$. So, $g$ has its
corresponding $h\in\iota(\sog(3))$ such that, $\Ad_g|_{\soa(3)}=\Ad_h|_{\soa(3)}$. Thus,
$\Ad_{gh^{-1}}|_{\soa(3)}={\rm Id}$, so that the element $g h^{-1}\in\sog(5)$ must satisfy
\[
g h^{-1} X = X g h^{-1},\quad \forall X\in \mathrm{Span}(E_1,E_2,E_3).
\]
Forcing $g h^{-1}$ to satisfy this condition on the basis $E_J$ for 
$J=1,2,3$, we find that $g h^{-1} = I$. Thus $g=h$ is in $G_\ten$ which
means that also $-g g^{-1}=-I$ is in $G_\ten$. But $-I\in\og(5)$ sends 
$\ten_{ijk}$ to $-\ten_{ijk}$, which gives the contradiction and finishes 
the proof.\\
\end{proof}

Thus we have the following proposition.
\begin{proposition}~\\
The stabiliser of tensor $\ten_{ijk}$ is the 
irreducible $\sog(3)$ included by $\iota$ in $\og(5)$. 
\end{proposition}

\subsection{The $\og(5)$ invariant characterisation of tensor $\ten$}

Since the stabiliser of $\ten_{ijk}$ is the irreducible $\sog(3)$, its
orbit under the $\og(5)$ action is a 7-dimensional 
homogeneous space $\og(5)/\iota(\sog(3))$. In this section we fully
characterise this orbit among all the orbits of $\og(5)$ action in 
$\bgs^3\bbR^5$. On doing this we view $\ten_{ijk}$ as 
a linear map 
$$\bbR^5\ni v\mapsto \ten_v\in\End(\bbR^5),\qquad\qquad
(\ten_v)_{ij} = \ten_{ijk} v_k.$$
Using this map we can rewrite the property iii) of Proposition
\ref{pr:ident} characterising $\ten_{ijk}$ to the equivalent form  
$$
\forall v\in\bbR^5\qquad \ten_v^2 v = g(v,v) v.$$
The importance of this reformulation is justified by the following theorem. 
\begin{theorem}
The $\og(5)$ orbit of tensor $\ten_{ijk}$ consists of all tensors
$\ft_{ijk}$ for which the associated linear map 
$$\bbR^5\ni v\mapsto \ft_v\in\End(\bbR^5),\qquad\qquad
(\ft_v)_{ij} = \ft_{ijk} v_k$$
satisfies the following three conditions
\begin{enumerate}
  \item it is totally symmetric, i.e. $g(u,\ft_v w) = g(w,\ft_v u) = g(u,\ft_w v)$,
  \item it is trace free {\rm tr}$(\ft_v)=0$,
  \item for any vector v $\in\bbR^5$
\begin{equation}
\label{tid}
\ft_v^2 v = g(v,v) v.
\end{equation}
\end{enumerate}
 \label{th:tchar}
\end{theorem}

\begin{remark}
The $\og(5)$ orbit of $\ten_{ijk}$, described invariantly in the above
theorem, consists of two disjoint $\sog(5)$ orbits: the orbit of
$\ten_{ijk}$ and the orbit of $-\ten_{ijk}$. Indeed, both tensors $\pm
\ten_{ijk}$ satisfy the three conditions characterising the $\og(5)$ orbit
and $\ten_{ijk}$ can not be sent to $-\ten_{ijk}$ via an element
$h\in\sog(5)$. Otherwise the element $-h$ preserves $\ten_{ijk}$
and as such belongs to $G_\ten$ which contradicts Lemma \ref{lem:ginso}.
\end{remark}


\begin{proof}[Proof of Theorem]
Let us consider tensor $\ten_{ijk}$ for which $\ten_{ijk}a^ia^ja^k$ has the
standard form (\ref{tensor}). Then its corresponding map $\ten_v$ in 
the $g$-orthonormal basis $\bas_i$ of (\ref{bass}), is
represented by the following matrices
\begin{equation}
  \label{tvmat}
\begin{gathered}
  \hspace{\stretch{1}} 
  \ten_{\bas_1}=\left(\begin{smallmatrix}
    \scm1&0&0&0&0\\
    0&1&0&0&0\\
    0&0&s&0&0\\
    0&0&0&1&0\\
    0&0&0&0& s\\
  \end{smallmatrix}\right) \qquad
  \ten_{\bas_2}=\left(\begin{smallmatrix}
    0&1&0&0&0\\
    1&0&0&0&0\\
    0&0&0&0&c\\
    0&0&0&0&0\\
    0&0&c&0&0\\
  \end{smallmatrix}\right) \hspace{\stretch{1}}  \\
\ten_{\bas_3}=\left(\begin{smallmatrix}
    0&0&s&0&0\\
    0&0&0&0&c\\
    s&0&0&\scm b&0\\
    0&0&\scm b&0&0\\
    0&c&0&0&0\\
  \end{smallmatrix}\right) \qquad
  \ten_{\bas_4}=\left(\begin{smallmatrix}
    0&0&0&1&0\\
    0&0&0&0&0\\
    0&0&\scm b&0&0\\
    1&0&0&0&0\\
    0&0&0&0&b\\
  \end{smallmatrix}\right) \qquad 
\ten_{\bas_5}=\left(\begin{smallmatrix}
    0&0&0&0&s\\
    0&0&c&0&0\\
    0&c&0&0&0\\
    0&0&0&0&b\\
    s&0&0&b&0\\
  \end{smallmatrix}\right),
\end{gathered}
\end{equation}
where $s=-\frac 12, \ b=c= \frac{\sqrt{3}}{2}$. The advantage
of introducing additional constant $b$ will be clear later in the
proof.\\

Now, let us take an arbitrary tensor $\ft_{ijk}$ satisfying the three 
assumptions of Theorem \ref{th:tchar}. The theorem will be proven if 
we manage to construct an orthonormal basis 
$(e_1,\ldots,e_5)$ in $\bbR^5$ in which the matrices $\ft_{e_j}$ take the same
form (\ref{tvmat}) as the matrices $\ten_{\bas_i}$.

\begin{lemma}
\label{lem:vw}
For any pair of orthogonal vectors $v,w$ the following identity holds
\[
g(v,v) w = 2 \ft_v^2 w + \ft_w \ft_v v.
\]
\end{lemma}
\begin{proof}[Proof of Lemma]
Applying (\ref{tid}) for the vector $v+ r w$ ($r\in \bbR$) we get
\[
r g(v,v) w + r^2 g(w,w) v = r \ft_v^2 w + r^2 \ft_w^2 v + r \ft_v \ft_w v + r^2 \ft_v \ft_w w + r \ft_w \ft_v v + r^2 \ft_w \ft_v w.
\]
The linear in $r$ term of this identity when compared with the
symmetry $\ft_w v = \ft_v w$ yields the thesis.\\
\end{proof}
\smallskip

The 5-th order homogeneous polynomial $\det(\ft_v)$ considered on the
unit sphere $\{v:g(v,v)=1\}$ satisfies 
$\det(\ft_{-v}) = - \det(\ft_{v})$. Thus, it can not have a fixed sign
everywhere on the sphere and there exists a unit vector $e_2$
such that $$\det(\ft_{e_2})=0.$$ 
Let 
\[
e_1 \colon = \ft_{e_2} e_2
\]
and let $e_4$ be the unit vector in the kernel of $\ft_{e_2}$:
\[
\ft_{e_2} e_4 =0.
\]
\begin{lemma}
The vectors $(e_1,e_2,e_4)$ are unit and pairwise orthogonal.
\end{lemma}
\begin{proof}
\begin{align*}
g(e_4,e_1) &= g(e_4,\ft_{e_2} e_2) = g(e_2,\ft_{e_2}e_4)=0,\\
g(e_4,e_2) &= g(e_4,\ft^2_{e_2} e_2) = g(e_2,\ft^2_{e_2}e_4)=0.
\end{align*}
Using Lemma \ref{lem:vw} for the unit orthogonal vectors $w=e_2$ and
$v=e_4$ we get $e_2 = \ft_{e_2} \ft_{e_4}e_4$ and so
\[
g(e_2,e_1)=g(e_2,\ft_{e_2} \ft_{e_2} \ft_{e_4}e_4)=g(\ft_{e_2}^2 e_2,\ft_{e_4}e_4)=0.
\]
Finally, the vector $e_1$ is unit:
\[
g(e_1,e_1)=g(\ft_{e_2} e_2,\ft_{e_2} e_2)=g(\ft_{e_2}^2 e_2,e_2)=1.
\]
\end{proof}
\medskip

The space $\mathrm{Span}(e_1,e_2,e_4)$ is $\ft_{e_2}$-invariant and $\ft_{e_2}$ restricted to this invariant space is trace-free; the same is true for the restriction of $\ft_{e_2}$ to the orthogonal complement $\mathrm{Span}(e_1,e_2,e_4)^\perp$. So, there exists a number $c\geq 0$ and a pair of unit vectors $(e_3,e_5)$ such that
\[
\ft_{e_2} e_3 = c\, e_5, \quad \ft_{e_2} e_5 = c\, e_3, \quad c\geq 0
\]
and the system $(e_1,e_2,e_3,e_4,e_5)$ is the orthonormal basis of $\bbR^5$. The matrix of $\ft_{e_2}$ in this basis has the form as in (\ref{tvmat}), but the constant $c$ is not fixed. 

Now, the use of the assumed properties of $(\ft_{ijk})$ 
and the successive application of 
Lemma \ref{lem:vw} proves that the matrices $\ft_{e_1}, \ldots,\ft_{e_5}$
have 
the form of (\ref{tvmat}) with the following 
restrictions to the constants $(b,c,s)$:
\[
s=-\frac 12\qquad c^2= \frac 34 \qquad b^2= c^2.
\]
If $b=-c$ then one can perform the following change of basis
\[
(e_1,e_2,e_3,e_4,e_5) \longmapsto (e_1,e_2,-e_3,-e_4,-e_5)
\]
resulting the change $b\mapsto (-b)$ in the matrices (\ref{tvmat}).

This finishes the proof of Theorem \ref{th:tchar}.\\
\end{proof}
\begin{corollary}
The tensor $\ten_{ijk}$ is fully determined by its properties listed in
Proposition \ref{pr:ident}. \label{cor:22}
\end{corollary}
\section{The $\sog(3)$ structure in $\bbR^5$ and the representations
  of $\sog(3)$}
The last corollary motivates the 
following definition.
\begin{definition}
An $\sog(3)$ structure on $\bbR^5$ is a pair $(g,\ten)$ where $g$ is a
Riemannian metric $g(A,A)=g_{ij}a^ia^j$ and $\ten$ is a ternary form
$\ten(A,A,A)=\ten_{ijk}a^ia^ja^k$ such that 
\begin{itemize}
\item[i)] $\ten_{ijk}=\ten_{(ijk)}$,
\vskip1mm
\item[ii)] $\ten_{ijj}=0$,
\vskip1mm
\item[iii)]$\ten_{jki}\ten_{lni}+\ten_{lji}\ten_{kni}+\ten_{kli}\ten_{jni}=g_{jk}g_{ln}+g_{lj}g_{kn}+g_{kl}g_{jn}.$ 
\end{itemize}
\end{definition} 

In this section we will use an $\sog(3)$ structure to define 
representations of $\sog(3)$ in 
$\bigotimes^2\bbR^5$. First, we recall the following well known theorem .
\begin{theorem}
\label{th:wigner}
All the irreducible finite-dimensional representations of $\sog(3)$
are odd dimensional. There is a unique irreducible representation of
$\sog(3)$ in space $\mathbb{R}^{2l+1}$ for each
$l\in\{0,1,2,3,...\}$. The tensor product
$\mathbb{R}^{2l_1+1}\otimes\mathbb{R}^{2l_2+1}$ decomposes onto the
$\sog(3)$-irreducible components according to the following Wigner formula
\begin{equation}
  \label{wigner}
  \mathbb{R}^{2l_1+1}\otimes\mathbb{R}^{2l_2+1} = \bigoplus_{l=|l_1-l_2|}^{|l_1+l_2|} \mathbb{R}^{2l+1}.
\end{equation}
\end{theorem}
\medskip

\noindent
The 5-dimensional irreducible representation $\rho$ of $\sog(3)$ with
the carrier space 
$$\bgw^1_5:=\bbR^5$$ was already considered in (\ref{so35}). To
find the projectors onto the irreducible components of the tensor
representations $\bigotimes^2\bbR^5$, $\bgw^2\bbR^5$ and
$\bgs^2\bbR^5$ we use the $\sog(3)$ structure $(g,\ten)$. Associated with
$\ten$ is the following endomorphism 
$$\hat{\ten}:{\textstyle \bigotimes}^2\bbR^5
\longrightarrow{\textstyle\bigotimes}^2\bbR^5,$$ 
$$W^{ik}\stackrel{\hat{\ten}}{\longmapsto} 4~\ten_{ijm}\ten_{klm} W^{jl},$$    
which preserves the decomposition
    $\bigotimes^2\bbR^5=\bgw^2\bbR^5\oplus\bgs^2\bbR^5$. 
Now, a simple algebra leads to the following proposition.
\begin{proposition} \label{pr:so3rep} $\bgt^2\bbR^5 = \bgw^2_3 \oplus \bgw^2_7 \oplus \bgs^2_1 \oplus \bgs^2_5 \oplus
  \bgs^2_9$, where
\begin{eqnarray*}
  &&\bgs^2_1       = \{~S\in\bgt^2\bbR^5~|~ \hat{\ten} (S)= 14 \cdot S~\} = \{S=\lambda
  \cdot g,\;\lambda\in\bbR~\}, \\
 && \bgw^2_3 = \{ ~F\in\bgt^2 \bbR^5~|~  \hat{\ten} (F) = 7 \cdot F~\}~ 
=~\soa(3)~=~\mathrm{Span}(E_1,E_2,E_3),\\
 && \bgs^2_5 = \{~ S\in\bgt^2\bbR^5 ~| ~ \hat{\ten}(S) = -3 \cdot S~\}, \\
 && \bgw^2_7 = \{ ~F\in\bgt^2\bbR^5 ~| ~  \hat{\ten}(F) = -8 \cdot F~\}=:\frak{n},\\
 &&\bgs^2_9 = \{~ S\in\bgt^2\bbR^5 ~| ~ \hat{\ten}(S) = 4\cdot S~\}.
\end{eqnarray*}
All the representations $\bgw^2_j\subset\bgw^2\bbR^5$ and $\bgs^2_k\subset\bgs^2\bbR^5$ are irreducible; the indices $j$ and $k$ denote their dimensions.
\end{proposition}
\begin{remark} Note that the tensor $\hat{\ten}$ defines a
nondegenerate $\sog(3)$ invariant scalar product 
$(F|F')=*(\hat{\ten}(F)\dz *F')$ of signature $(3,7)$ on the space of
2-forms 
\be
\bgw^2\mathbb{R}^5=\soa(5)=\soa(3)\oplus\frak{n}=\bgw^2_3\oplus\bgw^2_7.
\label{dec5}
\ee 
Although this scalar product
differs from the one associated with the 
Killing form $k(F,F')=-6*(F\dz *F')$, in both of them we 
have $\soa(3)\perp\frak{n}$.
\end{remark}
\begin{remark}
In agreement with the above notation we will denote the irreducible
representation $\bbR^5$ by $$\bgw^1_5 =\bbR^5 = \bgw^1 \bbR^5.$$
\end{remark}

\noindent
Using the $\sog(3)$ structure $(g,\ten)$ we can also build up an endomorphism  
$$ \check{\ten}:\bgs^2\bbR^5\longrightarrow\bgs^2\bbR^5$$
given by 
$$S^{kl}\stackrel{\check{\ten}}{\longmapsto} 4~\ten_{klm}\ten_{ijm} S^{ij}.$$
It is independent of $\hat{\ten}|_{\odot^2\bbR^5}$. 
Note that $\check{\ten}$ is a composition
$\check{\ten}=4\bar{\ten}\circ\grave{\ten}$ of two maps
$$
\bgs^2\bbR^5\stackrel{\grave{\ten}}{\longrightarrow}\bgw^1_5\stackrel{\bar{\ten}}{\longrightarrow}\bgs^2_5
$$
given by
\begin{equation}
\label{tltdef}
\grave{\ten}(S)_i=\ten_{ijk}S_{jk},\quad\quad\quad
\bar{\ten}(v)=\ten_{v}.
\end{equation}
We have 
$${\rm ker}(\grave{\ten})=\bgs^2_1\oplus\bgs^2_9,\quad\quad\quad {\rm
  im}(\grave{\ten})=\bgw^1_5.$$
Thus $\grave{\ten}$ restricted to $\bgs^2_5$ is an isomorphic intertwiner
  between the representations $\bgs^2_5$ and $\bgw^1_5$. Furthermore we have:
$$4\bar{\ten}\circ\grave{\ten}|_{\odot^2_5}=14\cdot {\rm id}.$$
Summarising we have the following proposition.
\begin{proposition}
The eigenvalues of $\check{\ten}$ on the representations
$\bgs^2_1\oplus\bgs^2_9$ and $\bgs^2_5$ are 0 and 14, respectively.
\end{proposition}

\section{The $\sog(3)$ structure on manifold}
\begin{definition}
An $\sog(3)$ structure on a 5-dimensional Riemannian manifold $(M,g)$
is a structure defined by means of a rank 3 tensor field $\ten$ for which
the associated linear map 
$${\rm T}M\ni v\mapsto \ten_v\in\End({\rm T}M),\qquad\qquad
(\ten_v)_{ij} = \ten_{ijk} v_k.$$
satisfies the following three conditions
\begin{enumerate}
  \item it is totally symmetric, i.e. 
$g(u,\ten_v w) = g(w,\ten_v u) = g(u,\ten_w v)$,
  \item it is trace free {\rm tr}$(\ten_v)=0$,
  \item for any vector field $v\in{\rm T}M$
$$
\ten_v^2 v = g(v,v) v.
$$
\end{enumerate}
\end{definition}

\begin{definition}
Two $\sog(3)$ structures $(M,g,\ten)$ and $(\bar{M},\bar{g},\bar{\ten})$
defined on
two respective 5-manifolds $M$ and $\bar{M}$ are (locally) \emph{equivalent}
iff there exists a (local) diffeomorphism $\phi:M\to \bar{M}$ such
that
$$\phi^*(\bar{g})=g\quad\quad{\rm and}\quad\quad \phi^*(\bar{\ten})=\ten.$$
If $\bar{M}=M$, $\bar{g}=g$, $\bar{\ten}=\ten$ the equivalence $\phi$ is
called a (local) \emph{symmetry} of $(M,g,\ten)$. The group of (local)
symmetries is called a \emph{symmetry group of} $(M,g,\ten)$.  
\end{definition}

\noindent
In view of Corollary \ref{cor:22}, Theorem \ref{th:tchar} and
Proposition \ref{pr:ident} tensor field $\ten$ reduces the structure
group of the bundle of orthonormal frames over $M$ to the irreducible
$\sog(3)$. Thus, locally, we can represent an $\sog(3)$ structure on
$M$ by a coframe 
\be
\theta=(\theta^i)=(\theta^1,\theta^2,\theta^3,\theta^4,\theta^5)\label{coframe}
\ee 
on
$M$, given up to the $\sog(3)$ transformation
\be
{\rm T}M\otimes\Om^1(M)\ni\theta\mapsto
\thet=\rho(h)\theta.\label{lifted}
\ee 
For such a class of coframes the Riemannian metric $g$ is
$$
g=\theta_1^2+\theta_2^2+\theta_3^2+\theta_4^2+\theta_5^2,
$$
and the tensor $\ten$, reducing the structure group from $\sog(5)$ to
$\sog(3)$, is 
\begin{equation}
\label{tform}
\ten=
\frac{1}{2}\theta_1\big(~6\theta_2^2+6\theta_4^2-2\theta_1^2-3\theta_3^2-3\theta_5^2~\big)+\frac{3\sqrt{3}}{2}\theta_4(\theta_5^2-\theta_3^2)+3\sqrt{3}\theta_2\theta_3\theta_5.
\end{equation}
\begin{definition}
\label{df:adapted}
An orthonormal coframe
$(\theta^1,\theta^2,\theta^3,\theta^4,\theta^5)$ in which
the tensor $\ten$ of an $\sog(3)$ structure $(M,g,\ten)$ is of the form
(\ref{tform}) is called a \emph{coframe adapted to $(M,g,\ten)$}, an
\emph{adapted coframe}, for short.
\end{definition}
\subsection{Topological obstruction}
The determination of topological obstructions for existence of an
irreducible $\sog(3)$ structure on a 5-dimensional manifold is presented in a
separate paper of one of us \cite{mbobi}. For the completeness of the
present paper we quote the result here. In the theorem below we denote 
by $p_j$ the $j$th Pontriagin class.
\begin{theorem}
\label{th:top}
Let $M$ be an orientable 5-dimensional manifold. There exists an 
irreducible $\sog(3)$ structure on $M$ iff $M$ admits the standard $\sog(3)$
structure (i.e. $TM$ splits on the rank 2 trivial bundle and a rank 3
complement) and  
\[ 
p_1(TM)= 5\; \tilde{p},\quad\text{where}\quad \tilde{p}\in H^4(M;\bbZ).
\]
\end{theorem}
\begin{remark}
The irreducible inclusion $\iota(\sog(3))\subset\sog(5)$ induces the
irreducible inclusion of $\tilde{\iota}({\bf Spin}(3))\subset{\bf Spin}(5)$. 
Assuming that $w_2(TM)=0$, the $\sog(5)$ structure on $M$ can be lifted to
the $\mathbf{Spin}(5)$ structure. It further may be reduced to the
$\tilde{\iota}(\mathbf{Spin}(3))$ structure on $M$, provided that $M$
admits an irreducible $\sog(3)$ structure.
\end{remark}

\subsection{$\soa(3)$ connection}
\label{sec:so3}
Given an $\sog(3)$ structure as above, we consider an $\soa(3)$
  connection on $M$ represented locally by means of an $\soa(3)$-valued 1-form
  $\Gamma$ given by 
\be
\Gamma=(\Gamma^i_{~j})=\gamma^1E_1+\gamma^2E_2+\gamma^3E_3,\label{charcon}
\ee
where $\gamma^1,\gamma^2,\gamma^3$ are 1-forms on $M$ and $E_I$ with 
  $I=1,2,3$ are given by (\ref{so3bas}).
This connection, having values in $\soa(3)\subset\soa(5)$, is
  necessarily metric. Via the Cartan structure equations,  
\be
\der\theta^i+\Gamma^i_{~j}\dz\theta^j=T^i
\label{chartor}
\ee
\be
\der\Gamma^i_{~j}+\Gamma^i_{~k}\dz\Gamma^k_{~j}=K^i_{~j},
\label{charcurv}
\ee
it defines the torsion 2-form 
$T^i$ and the $\soa(3)$-curvature 2-form $K^i_{~j}$. Using these forms we
define the torsion tensor $T^i_{~jk}\in(\bbR^5\otimes\bigwedge^2\bbR^5)$ 
and the $\soa(3)$-curvature tensor 
$r^I_{~jk}\in(\soa(3)\otimes\bigwedge^2\bbR^5)$, respectively, by 
$$T^i= \frac{1}{2}T^i_{~jk}\theta^j\dz\theta^k$$
and
\be
r^I=\der\gamma^I+\tfrac12\epsilon^I_{~JK}\gamma^J\dz\gamma^K=\frac{\sqrt{3}}{2}r^I_{~jk}\theta^j\dz\theta^k.\label{maler}
\ee
(Note that, $K=(K^i_{~j})=r^1 E_1+r^2 E_2+r^3 E_3.$) 
The connection
satisfies the first Bianchi identity
\be
K^i_{~j}\dz\theta^j=DT^i\label{bianchi1}
\ee
and the second Bianchi identity
\be
DK^i_{~j}=0,\label{bianchi2}
\ee
with the covariant differential defined by
$$DT^i=\der T^i+\Gamma^i_{~j}\dz T^j,\quad\quad\quad DK^i_{~j}=\der
K^i_{~j}+\Gamma^i_{~k}\dz K^k_{~j}-K^i_{~k}\dz\Gamma^k_{~j}.$$

Since the irreducible $\sog(3)$ was defined by the demand that it
preserves $g$ and $\ten$ we have the following proposition.
\begin{proposition}
Every $\soa(3)$ connection $\Gamma$ of (\ref{charcon}) is metric 
$$\stackrel{\Gamma}{\nabla}_v(g)\equiv 0$$
and preserves tensor $\ten$
$$\stackrel{\Gamma}{\nabla}_v(\ten)\equiv 0\quad\quad\quad
\forall v\in{\rm T}M.$$
\end{proposition}
\noindent

\subsection{$\sog(3)$ structures with vanishing torsion.}
\label{sec:zertor}
In this section we find all $\sog(3)$ structures $(M,g,\ten)$ which
admit $\soa(3)$ connections $\Gamma$ of (\ref{charcon}) with vanishing
torsion 
\be
T^i\equiv 0.\label{vantor}
\ee 
 
Assuming that $T^i$ is identically zero
and using the first Bianchi identity (\ref{bianchi1}) for $\Gamma$ we
easily find that a lot of 
components of the $\soa(3)$-curvature $r^I$ 
vanish. Explicitly, we find that in such a case the curvature forms
$(r^1,r^2,r^3)$ are expressible in terms of only one function
$r^1_{~15}$ and read
\be
r^1=r^1_{~15}\kappa^1,\quad\quad\quad\quad r^2=r^1_{~15}\kappa^2,
\quad\quad\quad\quad r^3=r^1_{~15}\kappa^3,\label{por}
\ee
where
\beq
&\kappa^1=\sqrt{3}\theta^1\dz\theta^5+\theta^2\dz\theta^3+\theta^4\dz\theta^5,\nonumber\\
&\kappa^2=\sqrt{3}\theta^1\dz\theta^3+\theta^2\dz\theta^5+\theta^3\dz\theta^4,
\nonumber\\
&\kappa^3=2\theta^2\dz\theta^4+\theta^3\dz\theta^5.\nonumber
\eeq
It further follows, that under the assumption of (\ref{vantor}), 
the second Bianchi identity (\ref{bianchi2}) implies 
that $$r^1_{~15}={\rm const}.$$
This means that $r^1_{~15}$ is a real parameter and that there is 
only 1-parameter family of $\sog(3)$ structures with vanishing
torsion. This family equips the principal fibre bundle $F(M)$ 
of $\sog(3)$ frames
$$\sog(3)\to F(M)\stackrel{\pi}{\to}M$$ over $M$ with an $\soa(3)$-connection 
\beq
\Gat&=&\rho(h)~\Gamma~\rho(h)^{-1}-\der\rho(h)~\rho(h)^{-1}\label{gh}\\
&=&\gat^1 E_1+\gat^2E_2+\gat^3 E_3.\nonumber
\eeq 
This, together with the lifted coframe 
\be
\thet=\rho(h)\theta\label{lcf}
\ee 
of
(\ref{lifted}), satisfies the following differential system
\beq
\der\thet^1&=&-\st\gat^1\dz\thet^5-\st\gat^2\dz\thet^3\nonumber\\
\der\thet^2&=&-\gat^1\dz\thet^3-\gat^2\dz\thet^5-2\gat^3\dz\thet^4\nonumber\\
\der\thet^3&=&\gat^1\dz\thet^2+\st\gat^2\dz\thet^1-\gat^2\dz\thet^4-\gat^3\dz\thet^5\nonumber\\
\der\thet^4&=&-\gat^1\dz\thet^5+\gat^2\dz\thet^3+2\gat^3\dz\thet^2
\label{stalakr}\\
\der\thet^5&=&\st\gat^1\dz\thet^1+\gat^1\dz\thet^4+\gat^2\dz\thet^2+\gat^3\dz\thet^3\nonumber\\
\der\gat^1&=&-\gat^2\dz\gat^3+r^1_{~15}\kat^1
\nonumber\\
\der\gat^2&=&-\gat^3\dz\gat^1+r^1_{~15}\kat^2
\nonumber\\
\der\gat^3&=&-\gat^1\dz\gat^2+r^1_{~15}\kat^3,
\nonumber
\eeq
where
\beq
\kat^1&=&\sqrt{3}\thet^1\dz\thet^5+\thet^2\dz\thet^3+\thet^4\dz\thet^5,\nonumber\\
\kat^2&=&\sqrt{3}\thet^1\dz\thet^3+\thet^2\dz\thet^5+\thet^3\dz\thet^4,\label{kappatil}
\\
\kat^3&=&2\thet^2\dz\thet^4+\thet^3\dz\thet^5.\nonumber
\eeq 
The {\it eight} linearly independent 1-forms
$(\thet^1,\thet^2,\thet^3,\thet^4,\thet^5,\gat^1,\gat^2,\gat^3)$
constitute a basis of 1-forms on the {\it eight} dimensional manifold
$F(M)$. Moreover, since equations (\ref{stalakr}) have 
only constant coefficients on
their right hand sides, the basis 
$(\thet^1,\thet^2,\thet^3,\thet^4,\thet^5,\gat^1,\gat^2,\gat^3)$ can
be identified with 
a basis of left invariant forms on a Lie group to which the bundle 
$F(M)$ is (locally) diffeomorphic. Thus we may identify $F(M)$ with a
local Lie group, the structure constants of which 
can be read off from the system (\ref{stalakr})-(\ref{kappatil}). We find that, depending
on the parameter $r^1_{~15}$, these structure constants correspond to 
\begin{itemize}
\item[i)] ${\bf SO}(3)\times_\rho\bbR^5$ group iff $r^1_{~15}=0$
\item[ii)] $\sug(3)$ group iff $r^1_{~15}>0$
\item[iii)] ${\bf SL}(3,\bbR)$ group iff $r^1_{~15}<0$. 
\end{itemize}

It further follows from the system (\ref{stalakr})-(\ref{kappatil})
that the tensors 
\be
\tilde{g}=\thet_1^2+\thet_2^2+\thet_3^2+\thet_4^2+\thet_5^2,\label{gt}
\ee
and  
\be
\tilde{\ten}=\frac{1}{2}\thet_1\big(~6\thet_2^2+6\thet_4^2-2\thet_1^2-3\thet_3^2-3\thet_5^2~\big)+\frac{3\sqrt{3}}{2}\thet_4(\thet_5^2-\thet_3^2)+3\sqrt{3}\thet_2\thet_3\thet_5\label{tt}
\ee
on $F(M)$ are preserved under the Lie transport 
along the fibres of $\sog(3)\to
F(M)\stackrel{\pi}{\to}M$. Moreover, these tensors are degenerate in precisely vertical
directions. Thus they descend 
to $M$ defining, respectively, $g$ and $\ten$, i.e. an $\sog(3)$ structure, there.
Locally, depending on the sign of $r^1_{~15}$, 
this structure is isomorphic to the homogeneous model $M_0$ 
in case i), the homogeneous model $M_+$ in case ii) and the
homogeneous model $M_-$ in case iii).  

\begin{theorem}
All $\sog(3)$ structures with vanishing torsion are locally isometric 
to one of the symmetric spaces $$M=G/{\bf SO}(3),$$ 
where $$G={\bf
  SO}(3)\times_\rho\mathbb{R}^5,\quad\quad \sug(3)\quad\quad {\rm or}\quad\quad {\bf SL}(3,\mathbb{R}).$$
The Riemannian metric $g$ and the tensor $\ten$ defining the $\sog(3)$
structure are obtained via (\ref{gt})-(\ref{tt})  
by means of  the left invariant forms
$(\thet^1,\thet^2,\thet^3,\thet^4,\thet^5,\gat^1,\gat^2,$ $\gat^3)$ on
$G$, which satisfy (\ref{stalakr})-(\ref{kappatil}). In all three cases the
metric $g$ is Einstein. It is flat in case of $G={\bf
  SO}(3)\times_\rho\mathbb{R}^5$. In the other two
cases the metric is not even conformally flat.
\end{theorem}

\noindent
\begin{proof} Only the last three sentences of the theorem remain to 
be proven. Since there is no torsion, 
the Levi-Civita connection for $g$, when
written in terms of the coframe 
$(\thet^1,\thet^2,\thet^3,\thet^4,\thet^5,\gat^1,\gat^2,\gat^3)$, is
simply $\Gat$ of (\ref{gh}). Then, the direct calculation shows
that the metric is Einstein with both the Ricci scalar and the Weyl
tensor being proportional, modulo a constant factor, 
to $r^1_{~15}$. \end{proof}

\begin{remark}
According to the last sentence of the theorem the spaces $M_\pm$ 
corresponding to nontrivial $\sog(3)$ structures without torsion are not of
constant curvature for the Levi-Civita connection of $g$.
\end{remark}

\begin{remark}
Note that 
\be
-\tilde{K}_0=\kat^I E_I \label{caco}
\ee 
is the curvature of the canonical
connection \cite{kono} on the symmetric space
$\sug(3)/\sog(3)$. Moreover, 
\label{rk:ccon}
the forms
$(\thet^1,\theta^2,\thet^3,\thet^4,\thet^5,\gat^1,\gat^2,\gat^3)$
define an absolute teleparalelism on $F(M)$. They can be collected to
an $\sua(3)$-valued matrix 
\be
\Gamma_{\rm Cartan}=\begin{pmatrix} 
0&\gat^3&\gat^2\\
&&\\
-\gat^3&0&\gat^1\\
&&\\
-\gat^2&-\gat^1&0\end{pmatrix}+i\begin{pmatrix} 
\frac{\thet^1}{\sqrt{3}}-\thet^4&\thet^2&\thet^3\\
&&\\
\thet^2&\frac{\thet^1}{\sqrt{3}}+\thet^4&\thet^5\\
&&\\
\thet^3&\thet^5&-2\frac{\thet^1}{\sqrt{3}}\end{pmatrix},\label{koncar}
\ee
which defines an $\sua(3)$-valued Cartan connection on the bundle 
$\sog(3)\to F(M)\to M$. The curvature of this connection 
\be
\Omega_{\rm Cartan}=\der\Gamma_{\rm Cartan}+\Gamma_{\rm
  Cartan}\dz\Gamma_{\rm Cartan} \label{cccn}
\ee
is
$$\Omega_{\rm Cartan}=(r^1_{~15}-1)\begin{pmatrix} 
0&\kat^3&\kat^2\\
&&\\
-\kat^3&0&\kat^1\\
&&\\
-\kat^2&-\kat^1&0\end{pmatrix}$$
and it vanishes iff the corresponding $\soa(3)$ connection $\Gamma$
has constant positive curvature determined by $r^1_{~15}=1$. 
\end{remark}
\begin{remark}
Remark \ref{rk:ccon} can be generalised leading to the description of  
$\sog(3)$ geometries with arbitrary $\soa(3)$ connection in terms of
an $\sua(3)$ Cartan connection on the fibre bundle $\sog(3)\to F(M)\to
M$. Indeed, given an $\sog(3)$
geometry with the adapted coframe
$(\theta^1,\theta^2,\theta^3,\theta^4,\theta^5)$ and the $\soa(3)$
connection $\Gamma$ we define the lifted coframe
$(\thet^1,\thet^2,\thet^3,\thet^4,\thet^5)$ via (\ref{lcf}) and the
1-forms $(\gat^1,\gat^2,\gat^3)$ via (\ref{gh}). Then, the $\sua(3)$-valued Cartan
connection on $F(M)$ is given by equation (\ref{koncar}).
The curvature (\ref{cccn}) of this connection 
satisfies the Bianchi identity
\be
{\rm D}\Omega_{\rm Cartan}=\der\Omega_{\rm
    Cartan}+\Gamma_{\rm Cartan}\dz\Omega_{\rm
    Cartan}-\Omega_{\rm Cartan}\dz\Gamma_{\rm Cartan}\equiv 0\label{bianchi}
\ee
and naturally splits onto the real and imaginary parts
$$
\Omega_{\rm Cartan}={\rm Re}(\Omega_{\rm Cartan})+i\sigma(\tilde{T})=\begin{pmatrix} 
0&\tilde{r}^3&\tilde{r}^2\\
&&\\
-\tilde{r}^3&0&\tilde{r}^1\\
&&\\
-\tilde{r}^2&-\tilde{r}^1&0\end{pmatrix}+i\begin{pmatrix} 
\frac{\Thet^1}{\sqrt{3}}-\Thet^4&\Thet^2&\Thet^3\\
&&\\
\Thet^2&\frac{\Thet^1}{\sqrt{3}}+\Thet^4&\Thet^5\\
&&\\
\Thet^3&\Thet^5&-2\frac{\Thet^1}{\sqrt{3}}\end{pmatrix}.
$$  
The imaginary part is simply the lift of the torsion $T$ of the
$\soa(3)$-connection $\Gamma$,
$$
\Thet=\rho(h)T.
$$
The real part can be collected to a $5\times 5$ matrix 
$$
\tilde{R}=\tilde{r}^1E_1+\tilde{r}^2E_2+\tilde{r}^3E_3.
$$
This satisfies 
\be
\tilde{R}=\tilde{K}-\tilde{K}_0,\quad\quad\quad\quad\tilde{K}=\rho(h)K\rho(h)^{-1},\label{rozn}
\ee
where $K$ is the $\soa(3)$ curvature of $\Gamma$ and  
$\tilde{K}_0$ is given by (\ref{caco}).
Thus, $\tilde{R}$ is the lift of the $\soa(3)$ curvature $K$ shifted by the
curvature $-\tilde{K}_0$ of the canonical connection on the
symmetric space $\sug(3)/\sog(3)$.
\end{remark}

\subsection{$\spin(3)$ connection}
The even Clifford algebra $Cl_0(5,0)$ has a 4-dimensional
faithful representation in
which the orthonormal vectors $(\bas_1,\bas_2,\bas_3,\bas_4,\bas_5)$ may be represented by 
\begin{equation}
  \label{spin}
\begin{gathered}
  \hspace{\stretch{1}} 
 \bas_1=\left(\begin{smallmatrix}
    0&0&1&0\\
    0&0&0&-1\\
    1&0&0&0\\
    0&-1&0&0\\
  \end{smallmatrix}\right) \qquad
 \bas_2=\left(\begin{smallmatrix}
    0&1&0&0\\
    1&0&0&0\\
    0&0&0&1\\
    0&0&1&0\\
  \end{smallmatrix}\right) \hspace{\stretch{1}}  \\
\bas_3=\left(\begin{smallmatrix}
    0&0&-i&0\\
    0&0&0&i\\
    i&0&0&0\\
    0&-i&0&0\\
  \end{smallmatrix}\right) \qquad
\bas_4=\left(\begin{smallmatrix}
    0&-i&0&0\\
    i&0&0&0\\
    0&0&0&-i\\
    0&0&i&0\\
  \end{smallmatrix}\right) \qquad 
\bas_5=\left(\begin{smallmatrix}
    1&0&0&0\\
    0&-1&0&0\\
    0&0&-1&0\\
    0&0&0&1\\
  \end{smallmatrix}\right),
\end{gathered}
\end{equation} 
One checks, by direct 
calculations, that
$$\bas_i^2=1,\quad\quad \bas_i\bas_j+\bas_j\bas_i=0,\quad\quad j\neq 
i=1,2,3,4,5.$$ 

Now, the double covering homomorphism ${\bf Spin}(5)\to \sog(5)$ 
induces the isomorphism of the Lie algebras $\spin(5)\to\soa(5)$. 
By means of this isomorphism an element $\bas_i\bas_j\in\spin(5)$,
$i<j$, is mapped to $(f_{ij})$ - a $5\times 5$
antisymmetric matrix having value 1 at its entry $f_{ij}$, value -1 at
$f_{ji}$ and value 0 in all the remaining 
entries.  
This implies that the basis of the Lie algebra $\spin(3)$
corresponding to the basis $(E_1,E_2,E_3)$ of the irreducible
$\soa(3)$ is 
$${\bf E}_1=\tfrac12(\sqrt{3}\bas_1\bas_5+\bas_2\bas_3+\bas_4\bas_5),\quad\quad
{\bf E}_2=\tfrac12(\sqrt{3}\bas_1\bas_3+\bas_2\bas_5+\bas_3\bas_4),$$
$$
{\bf E}_3=\tfrac12(2\bas_2\bas_4+\bas_3\bas_5).
$$
Explicitly:
$$
  {\bf E}_1=\tfrac12\left(\begin{smallmatrix}
    0&i&-\sqrt{3}&i\\
    i&0&-i&-\sqrt{3}\\
    \sqrt{3}&-i&0&-i\\
    i&\sqrt{3}&-i&0\\
  \end{smallmatrix}\right)\quad\quad
  {\bf E}_2=\tfrac12\left(\begin{smallmatrix}
    i\sqrt{3}&-1&0&-1\\
    1&i\sqrt{3}&-1&0\\
    0&1&-i\sqrt{3}&1\\
    1&0&-1&-i\sqrt{3}\\
  \end{smallmatrix}\right)
$$
\begin{equation}
{\bf E}_3=\tfrac12\left(\begin{matrix}
    2i&0&i&0\\
    0&-2i&0&i\\
    i&0&2i&0\\
    0&i&0&-2i\\
  \end{matrix}\right).\label{macdir}
\end{equation}

\noindent
Thus we have
$$
\spin(3)={\rm Span}({\bf E}_1,{\bf E}_2,{\bf E}_3)\quad\subset\quad \spin(5)={\rm
  Span}(\tfrac12\bas_i\bas_j, i<j=1,2,\dots,5).
$$
Now, given an $\sog(3)$ structure $(M,g,\ten)$ and an $\soa(3)$
connection $\Gamma=\gamma^1E_1+\gamma^2E_2+\gamma^3E_3$, we associate
with it a connection
\be\Gamma_{\spin}=\gamma^1{\bf E}_1+\gamma^2{\bf E}_2+\gamma^3{\bf
  E}_3\in \spin(3)\label{spicon}
\ee
which we call $\spin(3)$ connection. This connection will be used in
Section \ref{strominger} to define covariantly constant spinor
fields on $M$. 

\section{Characteristic connection}

Suppose now that we are given an $\sog(3)$ structure $(M,g,\ten)$ 
on a 5-dimensional manifold $M$. This defines the Levi-Civita
connection $\lc$ which, having values in $\soa(5)$, is an element
$\lc_{ijk}$ of 
$\soa(5)\otimes\bbR^5=\bgw^2\bbR^5\otimes\bbR^5$. In the following we
will be only interested in a subclass of $\sog(3)$ structures, which
we term {\it nearly integrable}.
\begin{definition}
An $\sog(3)$ structure $(M,g,\ten)$
is called \emph{nearly integrable} iff 
\be
(\stackrel{\scriptscriptstyle{LC}}{\nabla}_v \ten)(v,v,v)\equiv 0\label{ni}
\ee
for the Levi-Civita connection $\stackrel{\scriptscriptstyle{LC}}{\nabla}$.
\end{definition}
The condition (\ref{ni}), when written in an adapted coframe (\ref{coframe}), is
\be
\lc_{m(ji}~\ten_{kl)m}\equiv 0.\label{nearly}
\ee
This motivates an introduction of the map
$$\ten':\bgw^2\bbR^5\otimes\bbR^5\mapsto \bgs^4\bbR^5$$
such that
\beq
\ten'(\lc)_{ijkl}&=&12\lc_{m(ji}~\ten_{kl)m}\nonumber\\
&=&\lc_{mji}~\ten_{mkl}+\lc_{mki}~\ten_{jml}+\lc_{mli}~\ten_{jkm}\nonumber\\
&+&\lc_{mij}~\ten_{mkl}+\lc_{mkj}~\ten_{iml}+\lc_{mlj}~\ten_{ikm}\label{tprim}\\
&+&\lc_{mik}~\ten_{mjl}+\lc_{mjk}~\ten_{iml}+\lc_{mlk}~\ten_{ijm}\nonumber\\
&+&\lc_{mil}~\ten_{mjk}+\lc_{mjl}~\ten_{imk}+\lc_{mkl}~\ten_{ijm}.\nonumber
\eeq
We have the following proposition.
\begin{proposition}
An $\sog(3)$ structure $(M,g,\ten)$ is nearly integrable if and only if
its Levi-Civita connection $\lc\in\ker \ten'$. 
\end{proposition}
It is worthwhile to note that each of the last four raws of 
(\ref{tprim}) resembles the l.h.s. of equality (\ref{defso3}). Thus,  
$\soa(3)\otimes\bbR^5\subset\ker \ten'.$
Due to the first equality in (\ref{tprim}) we also have
$\bgw^3\bbR^5\subset\ker \ten'.$ It further follows that $\ker
\ten'=[\soa(3)\otimes\bbR^5]+\bgw^3\bbR^5.$ Now, introducing the map
$$\grave{\ten}:\ker \ten'\to\otimes^2\bbR^5$$ 
given by 
$$\grave{\ten}(\lc)_{il}=\ten_{ijk}\lc_{ljk}$$
and observing that $\ker\grave{\ten}=\bgw^3\bbR^5$ we get the following $\sog(3)$ invariant decomposition
$$\ker \ten'=[\soa(3)\otimes\bbR^5]\oplus\bgw^3\bbR^5.$$
This is the base for the following proposition.
\begin{proposition} \label{pr:charcon}
The Levi-Civita connection $\lc$ of a nearly integrable $\sog(3)$
structure $(M,g,\ten)$ uniquely decomposes onto
\be
\lc=\Gamma+\tfrac{1}{2}T,\label{rozk}
\ee
where 
$$\Gamma\in\soa(3)\otimes\bbR^5\quad\quad{\rm and}\quad\quad 
T\in\bgw^3\bbR^5=\ker\grave{\ten}.$$
\end{proposition}
The decomposition (\ref{rozk}) of the Levi-Civita connection $\lc$
of a nearly integrable $\sog(3)$ structure defines an $\soa(3)$
connection $\Gamma$. Rewriting 
the Cartan structure equation
$$\der \theta^i+\lc~^i_{~j}\wedge\theta^j=0$$
for $\lc$ into the
form 
$$\der
\theta^i+\Gamma^i_{~j}\wedge\theta^j=\tfrac{1}{2}T^i_{~jk}\theta^j\wedge\theta^k$$
enables us to interpret $T$ as the {\it totally skew symmetric
  torsion} of $\Gamma$.
\begin{definition}
An $\soa(3)$ connection $\Gamma$ of an $\sog(3)$ structure $(M,g,\ten)$
is called a \emph{characteristic connection} if its torsion $T_{ijk}$
is totally skew symmetric.
\end{definition}

The consideration of this section can be summarised in the following
theorem.

\begin{theorem}
Among all $\sog(3)$ structures only the nearly integrable ones 
admit characteristic connection $\Gamma$. Every nearly integrable
$\sog(3)$ structure defines $\Gamma$ uniquely.
\end{theorem}

\begin{remark}
Note, that out of {\it a priori} 50 independent components of the Levi-Civita
connection $\lc$, the nearly integrable condition (\ref{ni}) excludes
25. Thus, heuristically, the 
nearly integrable $\sog(3)$ structures constitute `a half' 
of all the possible $\sog(3)$ structures. 
\end{remark}

\begin{remark}\label{re:clastor}
Note, that given a nearly integrable $\sog(3)$ structure its totally
skew symmetric torsion $T_{ijk}$ defines the torsion 3-form
$$T=\tfrac{1}{6}T_{ijk}\theta^i\wedge\theta^j\wedge\theta^k.$$
Since $$\bgw^3\bbR^5=\bgw^2\bbR^5=\bgw^2_3\oplus\bgw^2_7$$ we have two 
kinds of skew symmetric torsions of `pure type' - those for which $T$
belongs to
the representation $\bgw^2_3$ and those whose $T$ is in $\bgw^2_7$.\\
Note that for an $\sog(3)$ structure with arbitrary $\soa(3)$
connection its torsion $T_{ijk}$ belongs to
$\bgw^2\bbR^5\otimes\bbR^5$. Thus, according to the discussion at the
beginning of this section, under the action of $\sog(3)$,
such $T_{ijk}$ satisfy 
$$T_{ijk}\in \bgw^2\bbR^5\otimes\bbR^5=\Big(
[\soa(3)\otimes\bbR^5]\oplus\bgw^3\bbR^5 \Big)\oplus \bbR^{25}.$$
Obviously, $\bbR^{25}$ further decomposes onto irreducibles: $\bbR^{25}=\bbR^5\oplus\bbR^9\oplus\bbR^{11}$.

\end{remark}

We close this section with the analysis of the $\sog(3)$ decomposition
of the curvature 
$$
K^i_{~j}=\tfrac{1}{2}K^i_{~jkl}\theta^k\dz\theta^l=\der\Gamma^i_{~j}+\Gamma^i_{~k}\dz\Gamma^k_{~j}
$$
of the characteristic connection $\Gamma$. Since
$K_{ijkl}\in\soa(3)\otimes\bgw^2\bbR^5$, this is given by the
following proposition. 
\begin{proposition}~
The projectors onto the irreducible components of the decomposition
\label{pr:curtordec}
  \begin{equation}
    \label{so3cur}
\soa(3)\otimes \bgw^2\bbR^5 \cong \bgs^2_1 \oplus \bgw^2_3 \oplus \bgw^2_7
\oplus \bgs^2_5 \oplus \bgs^2_9 \oplus \bgw^1_5,
\end{equation}
are:
\begin{eqnarray*}
&&K_{ijkl}\longmapsto  K_{[ijkl]} \in \bgw^4\bbR^5 = \bgw^1_5\\
&&K_{ijkl}\longmapsto  K_{ijil}=: k_{jl}\longmapsto 
  k_{[jl]}\in\bgw^2_3\oplus\bgw^2_7\\
&&K_{ijkl}\longmapsto  K_{ijil}= k_{jl}\longmapsto (
  k_{(jl)}-\frac{1}{5}k_{ii}~ g_{jl})\in\bgs^2_5 \oplus \bgs^2_9\\
&&K_{ijkl}\longmapsto  K_{ijil}= k_{jl}\longmapsto
 k_{ii}\in\bgs^2_1.
\end{eqnarray*}
\end{proposition}

\begin{remark}
Note that the curvature tensor decomposition
(\ref{so3cur}) is an analog, but not just the refinement, of the 
standard Riemann
tensor components. The $\soa(3)$-connection, we investigate, is
\emph{not} in general (compare Section \ref{sec:zertor} for the exception) the
torsion free connection and so the curvature does not have the usual Riemann tensor symmetries.
\end{remark}
\section{Homogeneous examples}
\label{sec:ex}

In the present section we look for examples of nearly integrable 
$\sog(3)$ structures admitting transitive symmetry groups.

Using the fact that the possible subgroups of $\sog(3)$ may have
dimensions 0,1,3 we get the following proposition.
\begin{proposition}
A transitive symmetry group $G$ of an
$\sog(3)$ structure may have the following dimension: $5,6$ or $8$.
\end{proposition}

\subsection{Examples with 8-dimensional symmetry group}

If the group of transitive symmetries $G$ is 8-dimensional, the 
$\sog(3)$ frame bundle $F(M)$ may be identified with
$G$. Then, the problem of finding all the examples with such group of
symmetries is equivalent to find those $G$s for which the basis of
left invariant forms
$(\thet^1,\thet^2,\thet^3,\thet^4,\thet^5,\gat^1,\gat^2,\gat^3)$ 
satisfy the pull-backed Cartan equations
(\ref{charcon})-(\ref{charcurv}) with the torsion coefficients
  $T_{ijk}$ and the curvature coefficients $r^I_{~jk}$ constant on
  $G$. This is a purely algebraic problem with the following solution.

\begin{proposition}
\label{pr:8ex}
There are only three different examples of nearly integrable 
$\sog(3)$ geometries with 8-dimensional symmetry group. These are the
torsion-free models:
\[
M_+=\sug(3)/\sog(3), \qquad M_0=(\sog(3)\times_\rho\mathbb{R}^5)/\sog(3),\qquad M_-=\mathbf{SL}(3,\bbR)/\sog(3).
\]
\end{proposition}

\subsection{Examples with 6-dimensional symmetry group}

To obtain all the examples with 6-dimensional transitive symmetry groups we
do as follows. We further reduce the lifted 
system (\ref{charcon})-(\ref{charcurv}) from the $\sog(3)$ frame bundle
$F(M)$ to a 6-dimensional group $G$ fibred over
$M$. We will identify $G$ with the transitive symmetry group of the
considered structure. Thus, $M$ will be a homogeneous  space $$M=G/H$$ with
$H$ -  a 1-dimensional subgroup of $G$. 

The reduction of the lifted
system (\ref{charcon})-(\ref{charcurv}) from $F(M)$ to $G$ implies
that on $G$, the two of the connection 1-forms
$(\tilde{\gamma}^1,\tilde{\gamma}^2,\tilde{\gamma}^3)$, say
$\tilde{\gamma}^1$ and $\tilde{\gamma}^2$, must be $\bbR$-linearly dependent on
the lift of the adapted coframe
$(\thet^1,\thet^2,\thet^3,\thet^4,\thet^5)$. 
Thus, in such case, the basis for 1-forms on $G$ is
$(\thet^1,\thet^2,\thet^3,\thet^4,\thet^5,\gat^3)$. It is subject to
the lift of the structure equations
(\ref{charcon})-(\ref{charcurv}). One of the integrability conditions
for these equations require that $\gat^1$ and $\gat^2$ must be of the
form 
\begin{equation}
  \label{6concon}
  \begin{aligned}
\gat^1 &= -b\, \thet^3 + a\, \thet^5,\\
    \gat^2 &= a\, \thet^3 + b\, \thet^5,
     \end{aligned} 
\end{equation}
where $a,b\in \bbR$. 

Due to the fact that all the coefficients in the pullback of the
Cartan structure equations (\ref{charcon})-(\ref{charcurv}) 
are constant on $G$, the closure of these equations implies the
following proposition.
\begin{proposition}
\label{pr:6ex}
All $\sog(3)$ nearly integrable structures with 6-dimensional 
symmetry group have (skew symmetric) torsion of the form
\[
T = t_1 \theta^1\wedge \theta^2\wedge \theta^4 + t_2 \theta^1\wedge\theta^3\wedge\theta^5. 
\]
There are three families of such geometries
\begin{enumerate}
\item\label{it:61} $b=t_1=t_2=0$, and $a$ arbitrary;
\item\label{it:62} $a=b=0$ and $t_1,\,t_2$ arbitrary;
\item\label{it:63} $a=0,\ b=\frac{t_1-2 t_2}{2\sqrt{3}}$ and $t_1,\,t_2$ arbitrary.
\end{enumerate}
\end{proposition}

Below we discuss all possibilities.

\noindent
{\bf The point \ref{it:61} of Proposition \ref{pr:6ex}.}
In this case the torsion is obviously zero and the $\soa(3)$ curvature form is
$$
  K= -a^2 \Big[\kappa^1\cdot E_1 +\kappa^2\cdot E_2 + \kappa^3\cdot E_3 \Big],$$
where $\kappa^1,\kappa^2,\kappa^3$ are given by (\ref{por}). 
Thus, in this case, we reconstruct two of the three torsion-free
  examples. For $a=0$ the respective $\sog(3)$ structure is equivalent
  to $M_0$. For $a\neq 0$ we reconstruct the structure 
$M_-=\mathbf{SL}(3,\bbR)/\sog(3)$. The latter case corresponds to 
the following 6-dimensional subgroup of $\slg(3,\bbR)$
\begin{equation*}
  G=\left\{M=\left( 
      \begin{smallmatrix}
        d&e&f\\
        g&h&k\\
        0&0&m\\
      \end{smallmatrix}
\right):\quad \det M =1 \right\}, \qquad H=\sog(2)=\left\{\left( 
      \begin{smallmatrix}
        \cos t&\sin t&0\\
        -\sin t& \cos t&0\\
        0&0&1\\
      \end{smallmatrix}
\right)\right\}.
\end{equation*}

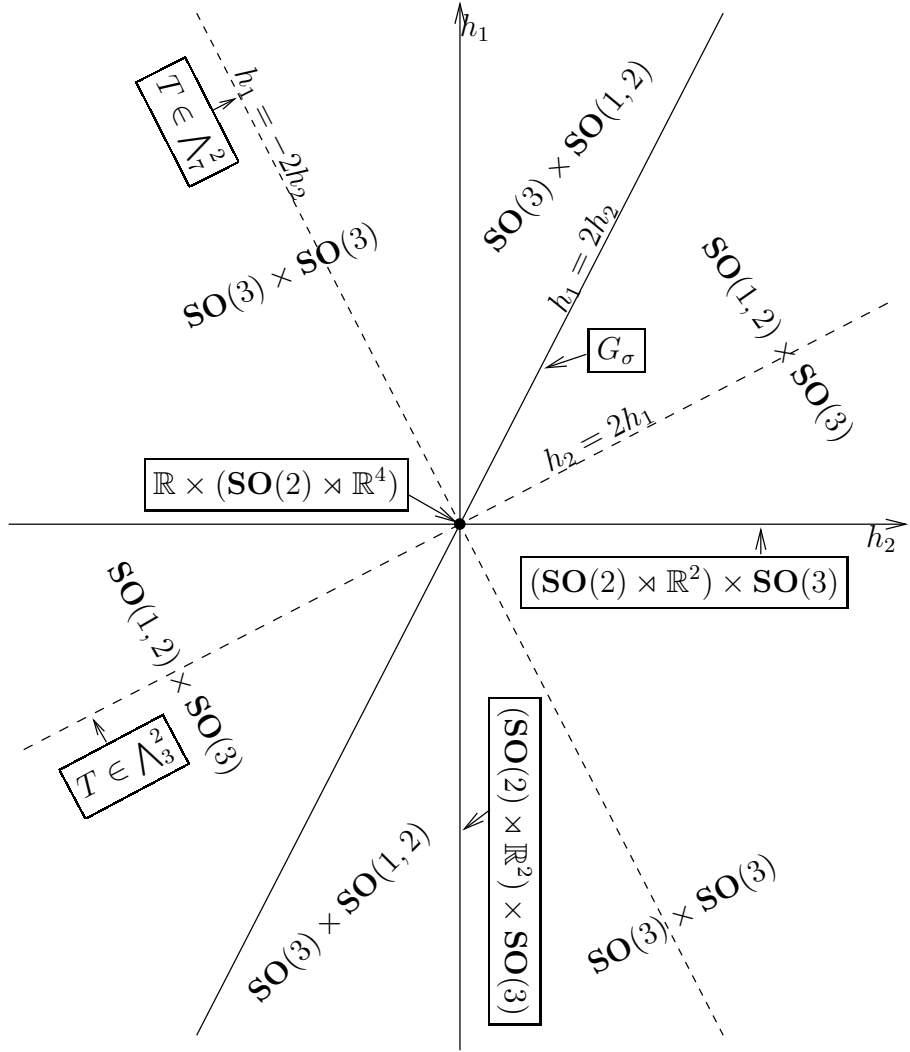
\begin{figure}[htbp]
  \begin{center}
    \input{hhex.pstex_t}
    \caption{Groups $G_{(t_1,t_2)}$ of $\sog(3)$ structures of
    Proposition \ref{pr:6ex} (\ref{it:62})}
    \label{fig:hhgroup}
  \end{center}
\end{figure}

\noindent
{\bf The point \ref{it:62} of Proposition \ref{pr:6ex}.}
In this case an invariant coframe $(\thet^1,\ldots,\thet^5,\gat^3)$ on
$G$ satisfies the following differential system:
\begin{align*}
  \dr \thet^1 &= t_1 \thet^2\wedge \thet^4 + t_2 \thet^3\wedge \thet^5\\
  \dr \thet^2 &= -t_1 \thet^1\wedge \thet^4 + 2 \thet^4\wedge \gat^3\\
  \dr \thet^3 &= -t_2 \thet^1\wedge \thet^5 + \thet^5\wedge \gat^3\\
  \dr \thet^4 &= t_1 \thet^1\wedge \thet^2 - 2 \thet^2\wedge \gat^3\\
  \dr \thet^5 &= t_2 \thet^1\wedge \thet^3 - \thet^3\wedge \gat^3\\
  \dr \gat^3  &= -
\frac{t_1 t_2}{2}(\thet^3\wedge\thet^5 + 2\thet^2\wedge\thet^4).
\end{align*}
The symmetry group $G=G_{(t_1,t_2)}$ depends on the torsion parameters $(t_1,
t_2)$. We depict the possible $G$s on the $(t_1,t_2)$-plane in Figure 
\ref{fig:hhgroup}.
\newcommand{\minus}{-}

Below we discuss each $G_{(t_1,t_2)}$ separately.

\begin{romlist}
\item $t_1\, t_2\,(t_1-2 t_2)\neq 0$. In this case $G$ is always of
  the form 
\begin{equation}
\label{g1g2}
G=G_1\times G_2,
\end{equation}
where $G_j$ is either $\sog(3)$ or $\sog(1,2)$ -- see Figure
\ref{fig:hhgroup}. There is a standard inclusion of $\sog(2)$ in both
of the above groups. The inclusion of $H=\sog(2)$ in the product $G$
is 
given by
$
\sog(2)\ni h \longmapsto (h^2,h)\in G_1\times G_2.
$
We consider the standard 3-dimensional representations of $\soa(1,2)$
and $\soa(3)$ so that the Maurer-Cartan form $\thet_{MC}$ on $G$ is given by
\begin{equation}
  \label{mcgen}
  \thet_{MC}=\begin{pmatrix}
  0& c \alt^1 + 2 \ett & \alt^2 &0&0&0\\
  -(c \alt^1 + 2 \ett) & 0 & \alt^4 &0&0&0\\
  \epsilon_1 \alt^2 & \epsilon_1 \alt^4&0 &0&0&0\\
  0&0&0& 0 & -2 c \alt^1 + \ett & \alt^3 \\
  0&0&0& -(-2 c \alt^1 + \ett) &0& \alt^5\\
  0&0&0& \epsilon_2 \alt^3&\epsilon_2 \alt^5&0\\
    \end{pmatrix},
\end{equation}
where
\begin{align*}
  c &= \tfrac1{\sqrt{5}},\\
  \epsilon_1 &= -\sgn[t_1 (t_1-2 t_2)],\\
  \epsilon_2 &= \sgn[t_2 (t_1-2 t_2)]
\end{align*}
and $(\alt^i,\ett)$ is a left invariant coframe on $G$.
We have the following relations between $(\alt^i,\ett)$ and the
canonical coframe $(\thet^i,\gat^3)$:
\begin{align}
\gat^3 &= \ett - \tfrac{2 t_1 + t_2 }{\sqrt{5}(t_1 -2 t_2)}\cdot \alt^1 &
\thet^1 &= \tfrac{\sqrt{5}}{t_1 -2 t_2}\cdot \alt^1\nonumber\\
\thet^2 &= \sqrt{\tfrac{\epsilon_1}{-t_1(t_1-2 t_2)}}\cdot\alt^2 &
\thet^4 &=\sqrt{\tfrac{\epsilon_1}{-t_1(t_1-2 t_2)}} \cdot\alt^4\label{altet}\\
\thet^3 &=\sqrt{\tfrac{2 \epsilon_2}{t_2(t_1-2 t_2)}}\cdot\alt^3 &
\thet^5 &=\sqrt{\tfrac{2 \epsilon_2}{t_2(t_1-2 t_2)}}\cdot\alt^5.\nonumber
\end{align} 

Now, we take $(\tilde{g},\tilde{\ten})$ in the canonical form (\ref{gt}),
(\ref{tt}). These descend to the $\sog(3)$ structure $(g,\ten)$ on
$M=G/H$ due to the isotropy invariance of $(\tilde{g},\tilde{\ten})$. 
The $G$-invariant $\soa(3)$ connection $\Gamma$ on $M$ has the following form
\[
 \Gamma = \Gamma_0 -\tfrac{1}{5}(2t_1+t_2)\theta^1\cdot E_3,
\]
where $\Gamma_0$ is the canonical connection on the reductive
homogeneous space $G/H$ -- see \cite{kono}.
\begin{remark}
It is worth to notice that on the line $t_2=-2 t_1$ the connection
$\Gamma$ coincides with the canonical connection $\Gamma_0$. The
example from this line corresponding to $(t_1,t_2)=(\frac15,-\frac25)$ 
is due to Th. Friedrich \cite{tfp}.
\end{remark}
In general, the torsion $T$ has components in the both  
possible irreducible $\sog(3)$ representations $\bgw^2_3$ and
$\bgw^2_7$ (see Remark \ref{re:clastor}). On the line $t_2=2t_1$ the torsion is of pure type $\bgw^2_3$; on the line $t_1=-2t_2$ it is of pure type $\bgw^2_7$ -- see the Figure \ref{fig:hhgroup}.

The $\soa(3)$ curvature is of the form
\[
K = - t_1 t_2 \kappa^3 \cdot E_3.
\]
It belongs to $\soa(3)\otimes\soa(3)$.  
If $t_1t_2\neq 0$ the curvature has
non-zero values in all of the components
$\bgs^2_1\oplus\bgs^2_5\oplus\bgw^1_5$ of the irreducible 
decomposition (\ref{so3cur}). 
\item $t_1=0,\ t_2\neq 0$. The group $G_1$ of the previous case 
contracts and the symmetry group becomes 
\[
G=(\sog(2)\rtimes\bbR^2)\times\sog(3).
\]
The inclusion of $H=\sog(2)$ in the product $G$ is given by
\[
\sog(2)\ni h \longmapsto (h^2,h).
\]
The Maurer-Cartan form on $G$ has the form (\ref{mcgen}) with $\epsilon_1=0$. The relations (\ref{altet}) remain valid after passing to the limit $\tfrac{\epsilon_1}{t_1}\to -\sgn[(t_1-2 t_2)] = \sgn t_2$
\begin{align*}
\gat^3 &= \ett + \tfrac{1}{2\sqrt{5}}\cdot \alt^1 &
\thet^1 &= -\tfrac{\sqrt{5}}{2 t_2}\cdot \alt^1\\
\thet^2 &= \tfrac{1}{\sqrt{2 |t_2|}}\cdot\alt^2 &
\thet^4 &=\tfrac{1}{\sqrt{2 |t_2|}}\cdot\alt^4\\
\thet^3 &=\tfrac{1}{|t_2|}\cdot\alt^3 &
\thet^5 &=\tfrac{1}{|t_2|}\cdot\alt^5.
\end{align*} 
These define an $\sog(3)$ structure on $M=G/H$ in an analogous way as
in the previous case. The torsion $T\neq 0$ is never of a pure
type and the $\soa(3)$ curvature $K\equiv 0$.

\item $t_2=0,\ t_1\neq 0$. This case is the same as the previous one. One has to put $\epsilon_2=0$ in (\ref{mcgen}) and $\tfrac{\epsilon_2}{t_2}\to \sgn t_1$ in (\ref{altet}). The statements about curvature and torsion are the same as in the previous point.
\item $t_1=0,\ t_2= 0$. In this case both the torsion and the
  $\soa(3)$ curvature vanish. Thus, this case corresponds to the flat
  model $M_0$. Hence the symmetry group $G$ is extendable to
  $\sog(3)\times_\rho \bbR^5$. For the purpose of the next point it is 
  useful to analyse $G$ more carefully. Let $\tau$ be the standard
  representation of $\sog(2)$ in $\bbR^2$. In conform with the 
Figure \ref{fig:hhgroup} we
  observe that  
\[
G=\bbR\times ( \sog(2)\rtimes\bbR^4),\quad\quad H=\sog(2), 
\]
where the semi-direct product is taken with respect to the
representation $\tau^2\oplus\tau$ of $\sog(2)$ on $\bbR^4$. The Maurer-Cartan form $\thet_{MC}$ on $G$ is
\begin{equation}
\label{mcflat}
  \thet_{MC}=\left(\begin{smallmatrix}
  0& 2 \ett & \alt^2 &0&0&0&0\\
  \scm 2 \ett & 0 & \alt^4 &0&0&0&0\\
  0& 0&0 &0&0&0&0\\
  0&0&0& 0 & \ett & \alt^3&0 \\
  0&0&0& \scm \ett&0& \alt^5&0\\
  0& 0&0 &0&0&0&0\\
  0&0&0&0 &0&0&\alt^1\\
    \end{smallmatrix}\right).
\end{equation}
The relation between $(\alt^j,\ett)$ and $(\thet^j,\gat^3)$ 
is $\thet^j=\alt^j$ and $\gat^3=\ett$.

\item \label{it:5} $t_1=2 t_2,\ t_2\neq 0$. In this case the group $G=G_\sigma$ has
the following abstract description. We present the Lie algebra of $G$
  as a central extension by $\bbR$ of a 5-dimensional algebra $\frak{l}$. Let us recall (see \cite{homalg}) that such extensions are classified by closed 2-forms $\sigma\in\bgw^2\frak{l}^*$.

Let $L=\sog(2)\rtimes \bbR^4$ with the representation
$\tau^2\oplus\tau$ of $\sog(2)$ as in the previous point; $\frak{l}$ is the Lie
algebra of $L$. We take the Maurer-Cartan forms $(\alt^2,\alt^3,\alt^4,\alt^5,\ett)$, defined in (\ref{mcflat}), as the basis of the left invariant forms on $L$. One can check that the following 2-form on $L$
\begin{equation}
  \label{sigma}
  \tilde{\sigma}=\alt^3\wedge\alt^5 + 2 \alt^2\wedge\alt^4, \qquad \sigma\colon = \tilde{\sigma}_e\in\bgw^2\frak{l}^*
\end{equation}
is closed. 

We define the Lie algebra $\frak{g}=\frak{g}_\sigma$ as a central extension of $\frak{l}$ by $\bbR$
\begin{equation}
\label{liext}
0\rightarrow \bbR \longrightarrow \frak{g} \xrightarrow{\pi} \frak{l}\rightarrow 0
\end{equation}
characterised by the element $\sigma$. Let $G=G_\sigma$ be a Lie group with Lie algebra $\frak{g}_\sigma$. We extend the basis of left invariant forms on $L$ to the (left invariant) basis $(\alt^1,\alt^2,\alt^3,\alt^4,\alt^5,\ett)$ on $G$. The differential $\dr \alt^1$ is (see \cite{homalg})
\[
\dr \alt^1 = \tilde{\sigma}.
\]
The exact sequence of Lie algebras (\ref{liext}) has a partial splitting $s:\soa(2)\hookrightarrow \frak{g}$ (i.e.\ the composition $\pi\circ s$ is the inclusion of $\soa(2)$ into $\frak{l}$) which defines the inclusion $H=\sog(2)\subset G$.

Finally, the relation between this basis $(\alt^j,\ett)$ and the canonical coframe $(\thet^j,\gat^3)$ is as follows
\[
\gat^3=\ett-t_2^2\cdot \alt^1,\quad \thet^1=t_2\cdot \alt^1,\quad \thet^2=\alt^2,\quad \thet^3=\alt^3,\quad \thet^4=\alt^4,\quad\thet^5=\alt^5.
\]
These define a nearly integrable $\sog(3)$ structure on $M=G/H$ as in
each of the previous cases. The torsion $T\neq 0$ is never of a pure type and the $\soa(3)$-curvature has the form 
\[
K = - 2 (t_2)^2 \kappa^3 \cdot E_3.
\]
\end{romlist}

\noindent
{\bf The point \ref{it:63} of Proposition \ref{pr:6ex}.}
We start with the observation that the line $t_1=2t_2$ on the 
$(t_1,t_2)$-plane in the present case and the line $t_1=2t_2$ of the previous case coincide 
(see Proposition \ref{pr:6ex}). Thus, in the entire analysis of this
case, we assume that $t_1\neq 2 t_2$.

We have the following differential system on $G$
\begin{align}
  \dr \thet^1 &= t_1 \thet^2\wedge\thet^4 + (t_1-t_2) \thet^3\wedge\thet^5\nonumber\\
  \dr \thet^2 &= -t_1 \thet^1\wedge\thet^4 + 2 \thet^4\wedge\gat^3\nonumber\\
  \dr \thet^3 &= -\tfrac12 t_1 \thet^1\wedge\thet^5 + \thet^5\wedge \gat^3 +\tfrac{t_1-2t_2}{2\sqrt{3}}\thet^2\wedge\thet^3+\tfrac{t_1-2t_2}{2\sqrt{3}}\thet^4\wedge\thet^5\nonumber\\
  \dr \thet^4 &= t_1 \thet^1\wedge\thet^2 - 2 \thet^2\wedge\gat^3\label{sy3}\\
  \dr \thet^5 &= \tfrac12 t_1 \thet^1\wedge\thet^3 - \thet^3\wedge\gat^3 -\tfrac{t_1-2t_2}{2\sqrt{3}}\thet^2\wedge\thet^5 -\tfrac{t_1-2t_2}{2\sqrt{3}}\thet^3\wedge\thet^4\nonumber\\
  \dr \gat^3  &= - \tfrac23 (t_1^2 -t_1 t_2+t_2^2)\thet^2\wedge\thet^4 - \tfrac12 t_1(t_1-t_2)\thet^3\wedge\thet^5.\nonumber
\end{align}
It follows, that off the line $t_1=2t_2$, independently of
$(t_1,t_2)$, the symmetry group $G=G_{\epsilon\sigma}$ is a central
extension of the group 
$$L=\slg(2,\bbR)\rtimes\bbR^2$$
by a 1-dimensional Lie group. 
$G_{\epsilon\sigma}$ is characterised by a closed 2-form
$\epsilon\sigma\in\bgw^2\frak{l}^*$, $\epsilon=\sgn |t_1-t_2|$.

It is convenient to choose the basis of left invariant forms
$(\alt^2,\alt^3,\alt^4,\alt^5,\ett)$ on $L$ so that the
Maurer-Cartan form $\thet_{MC}$ on $L$ reads
\[
\thet_{MC}=\begin{pmatrix}
-\alt^4 & \alt^2 + \ett & \alt^3\\
\alt^2 - \ett & \alt^4 & \alt^5\\
0&0&0\\
\end{pmatrix}.
\]
Obviously, we have $\sog(2)\subset\slg(2,\bbR)\subset L$.

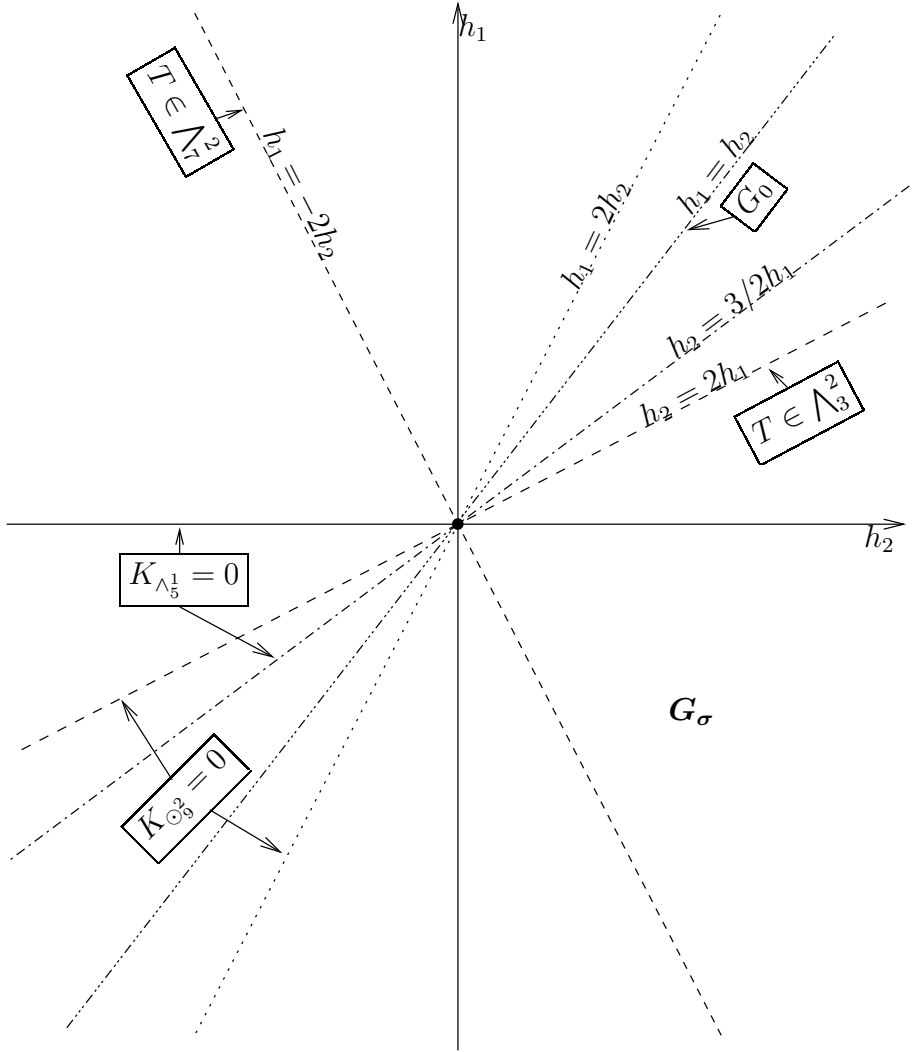
\begin{figure}[htbp]
  \begin{center}
    \input{hh2ex.pstex_t}
    \caption{Groups $G=G_{\epsilon\sigma}$ of $\sog(3)$ structures of
    Proposition \ref{pr:6ex} (\ref{it:63})}
    \label{fig:hh2g}
  \end{center}
\end{figure}

Now, the possible symmetry groups $G=G_{\epsilon\sigma}$,
$\epsilon=0,1$, are presented on Figure~\ref{fig:hh2g}.
Below, we discuss cases $\epsilon=1$ and $\epsilon=0$ separately.

\begin{romlist}
\item $\epsilon=1$. This case corresponds to $t_1\neq t_2$. Here, we
  observe that 
$$\tilde{\sigma}=\alt^3\wedge\alt^5,\quad\quad\quad \sigma\colon
  =\tilde{\sigma}_e\in\bgw^2\frak{l}^*$$
is closed on $L$. It is this form that defines the desired 
central extension of the Lie algebra $\frak{l}$ to the Lie algebra 
$\frg=\frg_\sigma$ of the symmetry group $G_\sigma$. Now, the forms 
$(\alt^2,\alt^3,\alt^4,\alt^5,\ett)$ extend to the left invariant
  forms on $G_\sigma$. Together with the form $\alt^1$ such that 
$\dr\alt^1=\tilde{\sigma}$
they define the left invariant coframe on $G_\sigma$. This coframe is
related to the canonical coframe $(\thet^i,\gat^3)$ of 
(\ref{sy3}) via
\begin{align}
\thet^1 &= -\tfrac{6 t_1}{(t_1-2t_2)^2}\cdot \ett + 
\tfrac{2(t_1-t_2)}{\epsilon}\cdot
\alt^1 & & \nonumber\\
\thet^2 &= \tfrac{2\sqrt{3}}{t_1-2t_2}\cdot \alt^2 &
\thet^4 &= \tfrac{2\sqrt{3}}{t_1-2t_2}\cdot \alt^4 \nonumber\\
\thet^3 &= \alt^3-\alt^5 & \thet^5 &= \alt^3+\alt^5\label{gat354}\\
\gat^3 &= \tfrac{(t_1-2t_2)^2+3t_1^2}{(t_1-2t_2)^2}\cdot\ett -
\tfrac{t_1(t_1-t_2)}{\epsilon}\cdot\alt^1. & &\nonumber
\end{align}
Now, in analogy to the case \ref{it:5} of page \pageref{it:5}, we use
the partial splitting $s:\soa(2)\to\frak{l}$, to recover the inclusion
  $H=\sog(2)\subset G_\sigma$. Then the $\sog(3)$ structure on $M=G/H$
is obtained via the standard procedure of taking
$(\tilde{g},\tilde{\ten})$ in the form (\ref{gt}), (\ref{tt}) and passing
to the quotient structure $(g,\ten)$. The $G$-invariant $\soa(3)$
connection on $M$ is given by
\[
\Gamma = \Gamma_0 -\tfrac{t_1-2t_2}{2\sqrt{3}}\theta^3\cdot E_1 +
\tfrac{t_1-2t_2}{2\sqrt{3}}\theta^5\cdot E_2 -
\tfrac{t_1}{2}\theta^1\cdot E_3, 
\]
where $\Gamma_0$ is the canonical connection on $G/H$.

As in the entire point \ref{it:62} of the present Proposition, the
torsion $T$ has the pure type $\bgw^2_3$ iff $t_2=2t_1$; it is of the
pure type $\bgw^2_7$ iff $t_1=-2t_2$; in all other cases it is not
of a pure type (see Figure \ref{fig:hh2g}).   

In contrast to the point \ref{it:62} of the present Proposition, the $\soa(3)$ curvature has the form
\begin{align*}
K &= \frac1{12} \left[\Big(\sqrt{3}t_1(t_1-2t_2)\, \theta^1\wedge\theta^5 - (t_1-2t_2)^2(\theta^2\wedge\theta^3+\theta^4\wedge\theta^5)\Big)\cdot E_1 \right.\\
 & + \Big(\sqrt{3}t_1(t_1-2t_2)\, \theta^1\wedge\theta^3 - (t_1-2t_2)^2(\theta^2\wedge\theta^5+\theta^3\wedge\theta^4)\Big)\cdot E_2  \\
    &\left. + \Big(-8(t_1^2 -t_1t_2 +t_2^2)\,\theta^2\wedge\theta^4 + (-7 t_1^2+10t_1t_2 -4t_2^2)\,\theta^3\wedge\theta^5 \Big)\cdot E_3\right],
\end{align*}
and (off the line $t_1=2t_2$) it is never of type
$\soa(3)\otimes\soa(3)$. In general, the curvature can assume values
in all of the components of the decomposition (\ref{so3cur}), but $\bgw^2_3$ and $\bgw^2_7$:
\[
K\in \bgs^2_1 \oplus \bgs^2_5 \oplus \bgs^2_9 \oplus \bgw^1_5.
\]
Independently of $(t_1,t_2)$ the curvature has always the $\bgs^2_1$
and $\bgs^2_5$ part; it is without the $\bgs^2_9$ component on the line
$t_2=2t_1$ and without the $\bgw^1_5$ component on lines $t_1=0$ and $3t_1=2t_2$ -- see Figure \ref{fig:hh2g}.

\item $\epsilon=0$. This corresponds to the line $t_1=t_2$. Now, all
  the formulas of the previous case remain valid, but the formulas for
  $\thet^1$ and $\gat^3$. To get correct expressions for them, one has
  to pass to the limit $\tfrac{t_1-t_2}{\epsilon}\to 1$ in
  (\ref{gat354}). 

It is worthwhile to note that the central extension
  $G_0$ is, in this case, trivial. Hence, the symmetry group is
  simply a product 
$$G_0=\bbR\times(\slg(2,\bbR)\rtimes\bbR^2)\quad\quad{\rm
    with}\quad\quad H=\sog(2)\subset\slg(2,\bbR).$$
\end{romlist}

\subsection{Examples with 5-dimensional symmetry group}
The first set of examples in this section is characterised by the 
requirement that a nearly integrable $\sog(3)$ geometry has flat characteristic
connection. The full list of such geometries is given in Section 
\ref{sec:nocur}. In Theorem \ref{th:nocur} we prove that flatness of
the characteristic connection implies that the corresponding nearly
integrable $\sog(3)$ geometry has at least 5-dimensional transitive symmetry
group. Inspection of the examples of Section \ref{sec:nocur} shows
that, in generic cases, their symmetry groups are strictly 5-dimensional.

Of course, examples with flat characteristic connections 
do not exhaust the list of all nearly integrable $\sog(3)$ structures 
with strictly 5-dimensional transitive symmetry group. We obtained
another two classes of examples assuming that, in addition to
the action of a 5-dimensional transitive symmetry group, the torsion of
characteristic connection is of pure type. The results are given in
respective Sections \ref{sec:tor23} and \ref{sec:tor27}. It is worth
noticing that it was possible to find \emph{all} structures with
5-dimensional transitive symmetry group and torsion in $\bgw^2_3$ (see
Theorem \ref{th:tor23}). In case of $\bgw^2_7$ type torsion we were
only able to find a 2-parameter family of examples.

\subsubsection{Vanishing curvature}
\label{sec:nocur}

\begin{theorem}
\label{th:nocur}
Let $(M,g,\ten)$ be a nearly integrable $\sog(3)$ structure with
vanishing curvature of its characteristic connection. Then $M$ has a
structure of a 5-dimensional Lie group $G$ and the $\sog(3)$-structure
is $G$-invariant.
\end{theorem}

\begin{proof}
Since the characteristic connection of an $\sog(3)$ structure is flat,
one can assume that the connection (locally) vanishes. Thus, in a
suitably chosen local coframe (\ref{coframe}) the first Cartan
structure equations are 
\begin{align}
\der \theta^1 &= t_1 \theta^2\dz\theta^3 + t_2 \theta^2\dz\theta^4 +
t_3 \theta^2\dz\theta^5 +t_4 \theta^3\dz\theta^4 +t_5
\theta^3\dz\theta^5 + t_6 \theta^4\dz\theta^5\nonumber\\
\der \theta^2 &= -t_1 \theta^1\dz\theta^3 - t_2 \theta^1\dz\theta^4 -
t_3 \theta^1\dz\theta^5 +t_7 \theta^3\dz\theta^4 +t_8
\theta^3\dz\theta^5 + t_9 \theta^4\dz\theta^5\nonumber\\
\der \theta^3 &= t_1 \theta^1\dz\theta^2 - t_4 \theta^1\dz\theta^4 -
t_5 \theta^1\dz\theta^5 -t_7 \theta^2\dz\theta^4 -t_8
\theta^2\dz\theta^5 + t_{10} \theta^4\dz\theta^5\label{ncsys}\\
\der \theta^4 &= t_2 \theta^1\dz\theta^2 + t_4 \theta^1\dz\theta^3 -
t_6 \theta^1\dz\theta^5 +t_7 \theta^2\dz\theta^3 -t_9
\theta^2\dz\theta^5 - t_{10} \theta^3\dz\theta^5\nonumber\\
\der \theta^5 &= t_3 \theta^1\dz\theta^2 + t_5 \theta^1\dz\theta^3 +
t_6 \theta^1\dz\theta^4 +t_8 \theta^2\dz\theta^3 +t_9
\theta^2\dz\theta^4 + t_{10} \theta^3\dz\theta^4.\nonumber
\end{align}
Here the functional coefficients $t_i, i=1,2,\ldots, 10$ are related
to the torsion 3-form $T$ via:
\begin{multline}
T= t_1 \theta^1\dz\theta^2\dz\theta^3+ t_2
  \theta^1\dz\theta^2\dz\theta^4+t_3
  \theta^1\dz\theta^2\dz\theta^5+ 
t_4 \theta^1\dz\theta^3\dz\theta^4+ \\
t_5 \theta^1\dz\theta^3\dz\theta^5+t_6
  \theta^1\dz\theta^4\dz\theta^5+ 
t_7 \theta^2\dz\theta^3\dz\theta^4+t_8
  \theta^2\dz\theta^3\dz\theta^5+t_9
  \theta^2\dz\theta^4\dz\theta^5+t_{10} \theta^3\dz\theta^4\dz\theta^5.
\end{multline}
Now, the Bianchi identities are equivalent to the following integrability
conditions of the system (\ref{ncsys}):
\begin{itemize}
\item[(a)] all the functions $t_i, i=1,2,\ldots, 10$ are {\it constants}
\item[(b)] they are subject to the following constraints
\begin{equation}
\label{dupablada}
\begin{split}
t_3 t_{10} + t_6 t_8 - t_5 t_9 &= 0\\
t_1 t_{10} + t_5 t_7 - t_4 t_8 &= 0\\
t_3 t_7  - t_2 t_8 + t_1 t_9 &= 0\\
t_2 t_{10} + t_6 t_7 - t_4 t_9 &= 0\\
t_{3} t_4 - t_2 t_5 + t_1 t_6 &= 0.
\end{split}
\end{equation}
\end{itemize}
The point (a) above proves the theorem, showing that $M$ can be
identified with the symmetry group $G$ which has $t_i$ as its
structure constants. \\
\end{proof}

Below we solve conditions 
(\ref{dupablada}) to fully characterise $G$ under the
genericity assumption  
$$t_{10}\neq 0.$$ 
If this is assumed the 
general solution of system (\ref{dupablada}) is 
\[
t_1 = \tfrac1{t_{10}} (t_4 t_8 - t_5 t_7),\quad t_2 = \tfrac1{t_{10}}
(t_4 t_9 - t_6 t_7),\quad t_3 = \tfrac1{t_{10}} (t_5 t_9 - t_6 t_8).
\] 
Now it is easy to see that the following linearly independent ($t_{10}\neq0!$) 1-forms 
\begin{align*}
\alpha^4 &= t_{10} \theta^1 - t_6 \theta^3 + t_5 \theta^4 -t_4
\theta^5\\
\alpha^5 &= t_{10} \theta^2 - t_9 \theta^3 + t_8 \theta^4 -t_7 \theta^5
\end{align*}
are closed. They can be further supplemented to a coframe
$(\alpha^1,\alpha^2,\alpha^3,\alpha^4,\alpha^5)$ on $M$ such that 
\begin{align*}
\der \alpha^1 &= \alpha^2\dz\alpha^3\\
\der \alpha^2 &= \alpha^3\dz\alpha^1\\
\der \alpha^3 &= \alpha^1\dz\alpha^2\\
\der \alpha^4 &= 0\\
\der \alpha^5 &= 0.
\end{align*}
This proves the following proposition.
\begin{proposition}
If the torsion coefficient $t_{10}\neq 0$, the symmetry group $G$ of
a nearly integrable $\sog(3)$-structure with flat characteristic
connection is isomorphic to $\sog(3)\times \bbR^2$.
\end{proposition}
\medskip

\subsubsection{Torsion in $\bgw^2_3$}
\label{sec:tor23}~\\

\noindent
In the following a parameter $\delta=0,1$ labels 5-dimensional Lie groups 
$G_\delta$. By definition $G_0=\sog(3)\times
Aff(1)$, the direct product of $\sog(3)$ and the affine group $Aff(1)$
in dimension 1. We characterise the group $G_1$ by specifying the structure
equations for a left invariant coframe on $G_1$. Thus, $G_1$ is such a
5-dimensional Lie group for which there exists a coframe 
$(\al^1,\al^2,\al^3,\al^4,\al^5)$ which satisfies the following equations:
\begin{align*}
  \der \alpha^1 &= 0\\
  \der \alpha^2 &= \alpha^1\wedge\alpha^2\\
  \der \alpha^3 &= -2\alpha^1\wedge\alpha^3\\
  \der \alpha^4 &= -\alpha^1\wedge\alpha^4 + \alpha^2\wedge\alpha^3\\
  \der \alpha^5 &= \alpha^2\wedge\alpha^4.
\end{align*}
This group has the Lie algebra $\frak{g_1}$ which is a central
extension $0\rightarrow \bbR \rightarrow \frak{g}_1 \rightarrow
\frak{h}\rightarrow 0$ of the 4-dimensional Lie algebra
\[
\frak{h}=\left\{\begin{pmatrix}x^1&x^3&x^4\\
0&-x^1&x^2\\0&0&0\\\end{pmatrix},~~x_1,x_2,x_3,x_4\in\bbR  \right\}
\]
by a real line $\bbR$. The extension is given by means of a closed $2$-form 
(see \cite{homalg}) $\sigma=\alpha^2\wedge\alpha^4$.

The following theorem is obtained by a successive application of the
Bianchi identities on the system (\ref{chartor})-(\ref{charcurv}) in
which the characteristic connection $\Gamma$ is supposed to have 
torsion in $\bgw^2_3$ and for 
which all the connection coefficients, the curvature coefficients and
the torsion coefficients are constants.
\begin{theorem}
\label{th:tor23}
Let $(M,g,\ten)$ be a nearly integrable $\sog(3)$ geometry admitting
a 5-dimensional transitive symmetry group $G$. Assume, in addition, that the
torsion of its characteristic connection is of pure type
$\bgw^2_3$. Then: 
\begin{itemize}
\item Modulo a constant $\sog(3)$ gauge transformation, it
is defined by means of the adapted coframe
$(\theta^1,\theta^2,\theta^3,\theta^4,\theta^5)$ satisfying the
following differential system:
\begin{equation*}
\begin{split}
  \der \theta^1 &= -\tfrac23\sqrt{3}\varrho\epsilon\Big(\theta^2\dz\theta^4 + (2-3\delta)
  \theta^3\dz\theta^5\Big)\\
  \der \theta^2 &= -2 \varrho \cos \varphi~ \theta^2\dz\theta^4\\
  \der \theta^3 &= -\varrho\cos\varphi~ \theta^2\dz\theta^5 +
  \sqrt{3}\varrho\epsilon(1-\delta)\theta^1\dz\theta^5
  + \varrho\epsilon\delta \theta^2\dz\theta^3 + \varrho(\delta\epsilon-\sin\varphi)\theta^4\dz\theta^5\\
\der \theta^4 &= -2 \varrho\sin\varphi~ \theta^2\dz\theta^4\\
\der \theta^5 &= \varrho\cos\varphi ~\theta^2\dz\theta^3 +
  \sqrt{3}\varrho\epsilon(\delta-1)\theta^1\dz\theta^3
  - \varrho\epsilon\delta \theta^2\dz\theta^5 -\varrho(\delta\epsilon+\sin\varphi) \theta^3\dz\theta^4,
\end{split}
\end{equation*}
with constant parameters $\varrho>0,\varphi\in [0,2\pi[,\epsilon=\pm1,\delta=0,1$.

\item $G\cong G_\delta$.
\item For all values of the parameters $\epsilon,\delta,\varrho,\varphi$
  the curvature of the characteristic connection is 
of type $\bgs^2_1\oplus\bgs^2_5\oplus\bgw^1_5$ with all the
  irreducible components non-zero.
\end{itemize}
\end{theorem}
\medskip

\subsubsection{Torsion in $\bgw^2_7$}~\\
\label{sec:tor27}

\noindent
It is easy to check that an adapted coframe
$(\theta^1,\theta^2,\theta^3,\theta^4,\theta^5)$ with differentials
given by:
\begin{equation*}
\begin{split}
\der\theta^1 &=0\\
  \der \theta^2 &= -\varrho\cos\varphi\theta^2\dz\theta^4\\
  \der \theta^3 &= \tfrac12\sqrt{3} \varrho \cos \varphi~
  \theta^1\dz\theta^3-\tfrac12\sqrt{3} \varrho \sin \varphi~
  \theta^1\dz\theta^5-\tfrac12 \varrho \sin \varphi~ \theta^2\dz\theta^3+
\tfrac12 \varrho \cos \varphi~ \theta^3\dz\theta^4\\
\der \theta^4 &= \varrho\sin\varphi~ \theta^2\dz\theta^4\\
  \der \theta^5 &= -\tfrac12\sqrt{3} \varrho \sin \varphi~
  \theta^1\dz\theta^3-\tfrac12\sqrt{3} \varrho \cos \varphi~
  \theta^1\dz\theta^5-\tfrac12 \varrho \sin \varphi~ \theta^2\dz\theta^5-
\tfrac12 \varrho \cos \varphi~ \theta^4\dz\theta^5,
\end{split}
\end{equation*}
where the parameters $\varrho>0,~\varphi\in [0,2\pi[$ are constants, 
defines a nearly integrable $\sog(3)$ geometry 
whose characteristic torsion has pure
type $\bgw^2_7$. Its symmetry group is transitive, strictly
5-dimensional and has the following Maurer-Cartan form 
\[
\theta_{MC}=\begin{pmatrix}
\al^4 & 0 &0 &\al^1\\
0&\al^5&0&\al^2\\
0&0&-(\al^4+\al^5)&\al^3\\
0&0&0&0
\end{pmatrix},
\]
where the forms $(\al^1,\al^2,\al^3,\al^4,\al^5)$
are related to the coframe 
$(\theta^1,\theta^2,\theta^3,\theta^4,\theta^5)$ via an appropriate 
$\varrho$-dependent ${\bf GL}(5,\bbR)$ transformation. 

It is worth noting that the curvature of the characteristic connection
in this 2-parameter family of examples is always of the type 
$\bgs^2_1\oplus\bgs^2_9$ with both the irreducible components non-zero.
\section{Ricci tensor and covariantly constant spinors}
\label{strominger}

\subsection{Ricci tensor}

We have the following proposition
\begin{proposition}
For every nearly integrable $\sog(3)$ structure $(M,g,\ten)$ the Ricci 
tensor $\riclc$ of the Levi-Civita connection $\lc$ is related to
the Ricci tensor $\ricga$ of the characteristic $\soa(3)$ 
connection $\Gamma$ via
$$\riclc_{ij}=\ricga_{ij}+\tfrac{1}{4}T_{ikl}T_{jkl}+
\tfrac{1}{2}(*\dr*T)_{ij}.$$
\end{proposition} 

\begin{corollary}
Given a nearly integrable $\sog(3)$ structure $(M,g,\ten)$ the 
following two conditions are equivalent.
\begin{itemize}
\item The codifferential of the torsion 3-form $T$
  vanishes.
\item The Ricci tensor $\ricga$ of the characteristic connection
  $\Gamma$ is symmetric.
\end{itemize}
\end{corollary}
Thus, for nearly integrable $\sog(3)$ structures we have
$$*\dr * T\equiv 0\quad\quad\iff\quad\quad
\ricga_{ij}\equiv \ricga_{ji}.$$

In the rest of this section we discuss the torsion/curvature
properties of the homogeneous examples of Section \ref{sec:ex}. 
It is interesting to note that all these examples  
satisfy $$*\dr * T\equiv 0.$$ Thus, the Ricci tensor
$\ricga$ is symmetric for them. In many cases both the Ricci tensors
$\ricga$ and $\riclc$ are diagonal\footnote{Note that the square of
  the matrix $E_3$ and
  its fourth power are diagonal matrices.}. 
More explicitly,
\begin{itemize}
\item in case (\ref{it:61}) of Proposition \ref{pr:6ex} we have:
$$\riclc=\ricga=-6a^2g,\quad\quad T\equiv 0$$
\item in case (\ref{it:62}) of Proposition \ref{pr:6ex} we have:
$$\riclc=\tfrac12
  (t_1^2+t_2^2)g+\tfrac{1}{24}(16t_1^2+12 t_1t_2-t_2^2)E_3^2+\tfrac{1}{24}(4t_1^2-t_2^2)E_3^4,$$
$$\ricga=\tfrac12 t_1t_2 E_3^2,$$$$\der T=-2t_1t_2\theta^2\dz\theta^3\dz\theta^4\dz\theta^5$$ 
\item  in case (\ref{it:63}) of Proposition \ref{pr:6ex} we have:
$$\riclc=
  (t_1^2-t_1t_2+\tfrac12 t_2^2)g+\tfrac{1}{24}(44t_1^2-58
  t_1t_2+27t_2^2)E_3^2+\tfrac{1}{24}(8t_1^2-10t_1t_2+3t_2^2)E_3^4,$$
$$\ricga=\tfrac12
  t_1(t_1-2t_2)g+\tfrac{1}{12}(14t_1^2-29
  t_1t_2+14t_2^2)E_3^2+\tfrac{1}{12}(t_1-2t_2)(2t_1-t_2)E_3^4,$$
$$\der T=-t_1^2\theta^2\dz\theta^3\dz\theta^4\dz\theta^5$$ 
\item for the examples of Theorem \ref{th:nocur} we have: 
$$\ricga\equiv 0,$$
$$\der T\equiv 0$$
and $\riclc$ has a rather complicated form depending on the torsion parameters
$t_a$, $a=1,2,\dots 10$; for some values of the parameters the
Levi-Civita Ricci tensor $\riclc$
may be diagonal, e.g.: if $t_a=0,\forall
a\neq 1$ then $$\riclc=\tfrac12 t_1^2\left(\begin{smallmatrix}
1 & 0 &0 &0&0\\
0&1&0&0&0\\
0&0&1&0&0\\
0&0&0&0&0\\
0&0&0&0&0\\
\end{smallmatrix}\right)$$
\item for the examples of Theorem \ref{th:tor23} we have:
$$\riclc=\varrho^2(\tfrac{10}{3}-2\delta) g+2\varrho^2 E_3^2,$$
$$\ricga=-2\varrho^2\delta g+\tfrac{4}{3}\varrho^2 E_3^2,$$
$$\der T=\tfrac43\varrho^2(3\delta-4)\theta^2\dz\theta^3\dz\theta^4\dz\theta^5$$
\item for the examples of Section \ref{sec:tor27} we have:
$$\riclc=-\tfrac{3\varrho^2}{2}\left(\begin{smallmatrix}
1 & 0 &0 &0&0\\
0&\sin^2(\varphi)&0&\tfrac12\sin(2\varphi)&0\\
0&0&0&0&0\\
0&\tfrac12\sin(2\varphi)&0&\cos^2(\varphi)&0\\
0&0&0&0&0\\
\end{smallmatrix}\right),\quad\ricga=-\tfrac{\varrho^2}{2}\left(\begin{smallmatrix}
3 & 0 &0 &0&0\\
0&2-\cos(2\varphi)&0&\sin(2\varphi)&0\\
0&0&1&0&0\\
0&\sin(2\varphi)&0&2+\cos(2\varphi)&0\\
0&0&0&0&1\\
\end{smallmatrix}\right),$$
$$\der T\equiv 0.$$

\end{itemize}

\subsection{Absence of covariantly constant spinors}

We now pass to the question if a manifold with an
$\sog(3)$ structure $(M,g,\ten)$ and an $\soa(3)$ connection $\Gamma$ may
admit a covariantly
constant spinor field. We look for $\Psi: M\to \bbC^4$ such that
\be
\der\Psi+\Gamma_\spin \Psi=0,\label{covspin}
\ee
where $\Gamma_\spin$ is a spin connection (\ref{spicon}) corresponding
to $\Gamma$.\\

We use the curvature 
$$\Omega_\spin=\der\Gamma_\spin+\Gamma_\spin\dz\Gamma_\spin$$ of
$\Gamma_\spin$. This curvature is expressible in terms of the
curvature $K=\tfrac{\sqrt{3}}{2}r^I_{~jk}\theta^j\theta^kE_I$ of $\Gamma$
and the (Dirac) matrices ${\bf E}_I$ of (\ref{macdir}). We have 
$$\Omega_\spin=\tfrac{\sqrt{3}}{2}r^I_{~jk}\theta^j\dz\theta^k{\bf
  E}_I.$$
It is easy to see that the integrability
conditions for the equations (\ref{covspin}) are   
$$\Omega_\spin\Psi=0.$$
These equations should be satisfied for each element of the basis of 2-forms
$\theta^i\dz\theta^k$. Thus, they are equivalent to 
$$W_{ij}\Psi=0,\quad\quad\forall i<j=1,2,3,4,5$$
where $W_{ij}$ is a 4x4 matrix
$$W_{ij}=r^I_{~ij}{\bf E}_I.$$
This shows that an existence of a non-zero solution for $\Psi$ gives a
severe restrictions on the curvature $\Omega_\spin$. In particular,
this implies that  
\be
\det(W_{ij})=\det(r^I_{~ij}{\bf
  E}_I)=0 \quad\quad \forall i<j=1,2,3,4,5.\label{iic}
\ee
But 
$$\det(W_{ij})=\tfrac{9}{16}\Big((r^1_{~ij})^2+(r^2_{~ij})^2+(r^3_{~ij})^2\Big)^2.$$
Thus, equations (\ref{iic}) are satisfied only if {\it all} the
  curvature coefficients $r^I_{~ij}$ are zero. In such case 
$\Omega_\spin=0$, which means 
that the corresponding $\soa(3)$ connection $\Gamma$ is {\it
  flat}. This proves the
  following proposition.
\begin{proposition}
Let $(M,g,t)$ be a 5-dimensional $\sog(3)$ geometry equipped with an
$\soa(3)$ connection $\Gamma$. Then $(M,g,\ten)$ admits a covariantly
constant spinor field with respect to the corresponding $\spin(3)$
connection $\Gamma_\spin$ if the connection $\Gamma$ is flat.
If this condition is satisfied then, locally, one has a 4-parameter
family of constant spinors.
\end{proposition}

\section{The twistor bundle $\twist$}
\label{sec:twistor}

It is remarkable that each 5 dimensional manifold $M$ with an $\sog(3)$
structure $(g,\ten)$ on it defines a natural $2$-sphere bundle 
$\sph^2\to\twist\to M$. This bundle, which via analogy with the twistor
theory, we call the \emph{twistor bundle}, will be defined by recalling
that at every point $x$ of $M$ we have a distinguished subspace
$(\bgw^2_3)_x$ of those
2-forms that span the irreducible $\soa(3)$. Considered point by point,
spaces $(\bgw^2_3)_x$ form a rank 3 vector bundle $\bgw^2_3M$ over $M$ with the
following basis of sections 
\begin{align*}
  \kappa^1 & = \sqrt{3}\theta^1\dz\theta^5+\theta^2\dz\theta^3+\theta^4\dz\theta^5,\\
  \kappa^2&=  \sqrt{3}\theta^1\dz\theta^3+\theta^2\dz\theta^5+\theta^3\dz\theta^4,\\
  \kappa^3 & = 2\theta^2\dz\theta^4+\theta^3\dz\theta^5.
\end{align*}
Here we have used the adapted coframe 
$(\theta^1,\theta^2\theta^3,\theta^4,\theta^5)$ for $(M,g,\ten)$. It is
also convenient to note that the forms $(\kappa^1,\kappa^2,\kappa^3)$
are related to the basis $(E_1,E_2,E_3)$
of the irreducible $\soa(3)\subset\soa(5)$ via 
$\kappa^I=\tfrac12 E_{Iij}\theta^i\wedge \theta^j,$ $I=1,2,3,$
see (\ref{so3bas}).
\begin{definition}
The \emph{twistor bundle} over a 5-dimensional manifold $M$ equipped with an 
$\sog(3)$ structure $(g,\ten)$ is the 2-sphere bundle 
$\sph^2\to\twist\xrightarrow{\pi} M$ defined by 
\be\label{twist}
\twist= \left\{ \omega\in\bgw^2_3M\ :\ *(\omega\dz*\omega)=5\  \right\}.
\ee
\end{definition}
\begin{remark}
The constant 5 in the above normalisation means that
$\omega\in\bgw^2_3M$ iff $\omega=b_1\kappa^1+b_2\,\kappa^2+
b_3\,\kappa^3$ where $b_1^2+b_2^2+b_3^2=1$.
\end{remark}

Consider the complexification ${\rm T}^\bbC M$ of the tangent bundle of
$(M,g,\ten)$ and denote by the same letters the complexifications of the
tensors $g$ and $\ten$. At every point $x\in M$ consider the space  
$$N_x=\{n\in {\rm T}_x^\bbC M\ : \ \ten(n,n,\cdot)\equiv 0\}$$
of vectors, which are {\it null} with respect to the complexified $\ten$. Given
any complexified vector $0\neq v\in {\rm T}^\bbC_xM$ we define 
$${\rm dir}(v)=\{\lambda v\in {\rm T}^\bbC_x M ~:~ \lambda\in \bbC\}.$$
We have the following proposition
\begin{proposition}
The space of null directions 
$$\bbP N_x=\{{\rm dir}(n) \ : \ n\in N_x\}$$
is a disjoint sum of two connected components
\[
\bbP N_x = \bbP N_x^+ \sqcup \bbP N_x^-,\quad\quad\quad\quad \bbP N_x^- =
\overline{\bbP N_x^+}.
\]
Each of them is naturally diffeomorphic to the fibre
$\twist_x=\pi^{-1}(x)=\sph^2$ of the twistor bundle $\twist$.
\end{proposition}
\begin{proof}
Consider a 2-form $\omega\in\twist_x$. In the adapted coframe
$(\theta^1,\theta^2,\theta^3,\theta^4,\theta^5)$ it reads
$\omega=\tfrac{1}{2}\omega_{ij}\theta^i\dz\theta^j$. It defines a linear map 
\[
\mathrm{T}_x^\bbC M \ni v_i \longmapsto (\omega v)_j = \omega_{ji}v_i
\in \mathrm{T}_x^\bbC M.
\]
It is easy to see that the eigenvalues of this endomorphism are
$\{0,\pm i,\pm 2 i\}$. The corresponding $\pm 2i$ eigenspaces are {\it
  null} with respect to $\ten$ due to the following argument. 
The
$\sog(3)$ invariance of the tensor $\ten$, see (\ref{defso3}), when
applied to form $\omega$ and a vector $n$ belonging to the $\pm2 i$
eigenspaces of $\omega_{ij}$ reads
\[
0=\ten(\omega n,n,\cdot) + \ten(n,\omega n,\cdot) + \ten(n,n,\omega \cdot)= 4i \ten(n,n,\cdot) +
\ten(n,n,\omega \cdot).
\]
Now, if $v$ belongs to any eigenspace of $\omega_{ij}$ the implication
of this equality is $\ten(n,n,v)=0$, 
which means that $\ten(n,n,\cdot)\equiv 0$. 

Thus, the following map 
\[
\twist_x \ni \omega\longmapsto \ker (\omega \mp 2 i) \in \bbP N_x^\pm.
\]
is well defined. It further follows that it provides the desired
diffeomorphism between $\twist_x$ and $\bbP N_x^\pm$.\\
\end{proof}
 
Now we define the 2-sphere {\it bundle of null directions} for $\ten$ to be 
$$ \bbP N=\bigcup_{x\in M} \bbP N_x^+$$
and, as a corollary to the above proposition, we get:
\begin{proposition}
There exists a natural bundle isomorphism between the bundle $\bbP N$
of null directions for $\ten$ and the twistor bundle $\twist$.
\end{proposition}

\begin{remark}
The above proposition enables one to view the twistor bundle $\twist$ 
as an analog of the twistor bundles of 4-dimensional
(pseudo)Riemannian geometries (see e.g. \cite{nur}). Historically, the 
first such bundle -
Penrose's bundle of light rays over the Minkowski space-time
\cite{pen} -  is a
2-sphere bundle of null directions. It proved to be very useful in
General Relativity Theory, especially in the case of complexified
Minkowski space-time and its curved generalisations. Motivated by the utility
of Penrose's bundle of light rays Atiyah, Hitchin and Singer
\cite{atiyah} considered the 2-sphere bundle of complexified null
2-planes over a 4-dimensional Riemannian manifold. This bundle,
which they identified with the bundle of almost hermitian structures
over the 4-manifold, they termed the {\it twistor bundle}. Later,
mathematicians generalised the notion of {\it twistor bundle} in many
directions, so that the relation between null directions and todays
twistors is weaker and weaker. We find a particularly remarkable the
fact that the 5-dimensional geometries considered in the present paper
lead to twistor bundle $\twist$ whose relation to null directions is
very apparent.
\end{remark}

\subsection{Elements of geometry of $\twist$}
Now, we consider an arbitrary $\sog(3)$ structure $(M,g,\ten)$ equipped with 
an $\soa(3)$ connection $\Gamma$ (we do not assume that $\Gamma$
is the characteristic connection). These data induce interesting
geometrical structures on the twistor bundle $\twist$. 
The rest of this section is devoted to their brief description. 
\begin{enumerate}
\item The connection $\Gamma$ splits the tangent space $\mathrm{T}\twist$ into horizontal and vertical parts:
\[
\mathrm{T}\twist = \calH\oplus \calV.
\]
The fibre of the twistor bundle $\twist_x$ is naturally embedded 
in the vector space $(\bgw^2_3)_x$. It is a unit sphere $\sph^2$ 
with respect to the natural scalar product $\Sigma$ on two-forms,
which explicitly reads
$$\Sigma(\sigma_1,\sigma_2)=\tfrac15*(\sigma_1\dz*\sigma_2),
\quad\quad \forall \sigma_1,\sigma_2\in(\bgw^2_3)_x.$$ 
Thus the vertical tangent space 
$\calV_\omega$ at a point $\omega\in\twist_x$ may be identified with 
the orthogonal complement of $\omega$ with respect to $\Sigma$.
Hence 
\[
\calV_\omega = \{\sigma \in(\bgw^2_3)_x\ :\quad \Sigma(\sigma,\omega) =0\}. \]
\item There is a natural Riemannian metric $\tilde{g}$ on
  $\twist$. This metric is given by
  $\tilde{g} = \Sigma_2\oplus \pi^*g$, where $\Sigma_2$ is the natural
  scalar product induced on the fibre by $\Sigma$.
\item There is a natural complex structure $J$ on the fibre
  $\twist_x$, given by 
\[
J_\omega (\sigma) = [\omega,\sigma]\quad\quad \sigma\in \calV_\omega\subset (\bgw^2_3)_x.
\]
Here, we view the forms $\omega$ and $\sigma$ as elements of the Lie
algebra $\soa(3)\cong(\bgw^2_3)_x$, so that $[\cdot,\cdot]$ is the Lie 
bracket in $\soa(3)$. Obviously, $J$ is
compatible with the metric $\Sigma_2$. Now, the metric $\Sigma_2$
together with orientation given by $J$ determine the volume 2-form
$\eta_2$ on the fibre.
\item There is a tautological horizontal 2-form $\omega$ on $\twist$.
\item $\twist$ is equipped with the horizontal vector field $u$
  given by 
\[
\tilde{g}(u)=\tfrac14\tilde{*}(\eta_2\dz\omega\dz\omega)
\]
where $\tilde{*}$ is the Hodge star operation on $(\twist,\tilde{g})$. The
vector field $u$ is unital: $\tilde{g}(u,u)=1$. We denote the
$\tilde{g}$-orthogonal complement of $u$ in $\calH$ by $\mathcal{H}^u$.
\item At every point $x\in\twist$ the metric $\tilde{g}$ descends to
  the 4-dimensional, naturally oriented, vector space
  $\calH^u_x$. Thus, in  $\calH^u_x$, the Hodge star
  operator is well defined. By using it we decompose the
  restriction of the tautological 2-form $\omega_{|_{\calH^u}}$ into
  the 
self-daul and anti-self-dual parts 
$$
\omega_{|_{\calH^u}}=\omega_+ + \omega_-.
$$
The forms $\omega_\pm$ define the pair of $\pi^*g$-compatible complex
 structures $J_\pm$ on $\calH^u$
\[
\pi^*g(J_\pm v_1, v_2)~=~\frac{2}{2\pm1}~\omega_\pm(v_1,v_2),\quad\quad
v_1,v_2\in \Gamma(\calH^u).
\]
These two complex structures commute: 
\[
[J_+,J_-]=0.
\]
\end{enumerate}

\subsection{Almost $CR$-structures on $\twist$ and their integrability
conditions}

We recall that an odd-dimensional real manifold $P$ is equipped with 
an {\it almost $CR$-structure} if there exists on $P$ a distinguished 
codimension one distribution
$\mathcal{N}$ endowed with an almost complex structure $\mathcal{J}$ 
(see e.g. \cite{nt}). The
$\pm i$ eigenspaces of $\mathcal{J}$ define the split 
$$\bbC\otimes\mathcal{N}=\mathcal{N}^{(1,0)}\oplus\mathcal{N}^{(0,1)}.$$
An almost CR-structure $(\mathcal{N},\mathcal{J})$ on $P$ is called an
{\it integrable} CR-structure iff the following integrability
conditions are satisfied
$$
[\mathcal{N}^{(1,0)},\mathcal{N}^{(1,0)}]\subset \mathcal{N}^{(1,0)}.
$$

The twistor bundle $\twist$ is naturally equipped with
{\it four} almost CR-structures. They are genuinely distinct i.e. not related
by the conjugacy operation. One obtains these
structures by defining the distribution $\mathcal{N}$ to be  
$\mathcal{N}_\twist=u^\perp,$
the orthogonal complement of the unit vector $u$ with respect to the
metric $\tilde{g}$ on $\twist$. Since 
$\mathcal{N}_\twist=\mathcal{V}\oplus\mathcal{H}^u$, then the four 
almost complex structures on $\mathcal{N}_\twist$ may be defined by
$$\mathcal{J}=J\oplus\epsilon J_\pm,\quad\quad\quad
\epsilon=1\ {\rm or}\ -1.$$
Thus we have four natural almost CR-structures on $\twist$ defined by
means of four $\mathcal{J}$s on $\mathcal{N}_\twist$. Among them the
most interesting is 
$$
(\mathcal{N}_\twist,\mathcal{J}_0), \quad {\rm where}\quad \mathcal{J}_0=J\oplus J_+.
$$
This structure is the only one among
$(\mathcal{N}_\twist,\mathcal{J})$ that may be integrable. 
More specifically, we have the following theorem. 
\begin{theorem}
\label{th:crint}~
\begin{enumerate}
\item Among the four natural almost CR-structures 
$(\mathcal{N}_\twist,J\oplus\epsilon J_\pm)$ on $\twist$, the only one 
that may be integrable is $(\mathcal{N}_\twist,\mathcal{J}_0)$. 
\item Let $(M,g,\ten)$ be a nearly integrable $\sog(3)$ structure and let 
  $(\mathcal{N}_\twist,\mathcal{J}_0)$ be the almost CR-structure on
  $\twist$ induced by the characteristic
  connection of $(M,g,\ten)$. This CR-structure is integrable if and only if
\[
K_{\odot^2_9} \equiv 0, \quad \text{and} \quad T\in\bgw^2_3.
\]
\end{enumerate}
\end{theorem}
\newcommand{\zb}{\overline{z}}
\begin{proof}[Sketch of the proof]
We start by choosing an $\sog(3)$ adapted
coframe 
$(\theta^1,\theta^2,\theta^3,\theta^4,\theta^5)$ on $U\subset M$. 
We parametrise $\twist|_{U}$ by $U\times \overline{\bbC}$, so
that the tautological 2-form $\omega$ reads
\begin{equation}
  \label{taut}
  \omega = \tfrac{z+\zb}{1+|z|^2}\kappa_1 + \tfrac{i(\zb-z)}{1+|z|^2}\kappa_2 + \tfrac{1-|z|^2}{1+|z|^2}\kappa_3,\qquad z\in \overline{\bbC}.
\end{equation}
The horizontal-vertical splitting of the tangent bundle $T\twist$ with respect to an $\soa(3)$-connection $\Gamma=\gamma_1 E_1 + \gamma_2 E_2 + \gamma_3 E_3$ is given by the following complex valued 1-form
\[
\tilde{h}=\tfrac{1}{1+|z|^2}\Big(\dr z + \tfrac{1-z^2}{2 i} \gamma_1 + \tfrac{1+z^2}{2} \gamma_2 + iz \; \gamma_3 \Big).
\]
The horizontal subspace $\mathcal{H}\subset T\twist$ is the kernel of 
$\tilde{h}$.

The 1-form $\tilde{u}=\tilde{g}(u)$ - the $\tilde{g}$-dual to the
unit horizontal vector field $u$ - is given by
\[
\tilde{u} = -\tfrac{1-4|z|^2+|z|^4}{(1+|z|^2)^2} \theta^1 +\tfrac{i\sqrt{3}(z-\zb)(z+\zb)}{(1+|z|^2)^2} \theta^2 - \tfrac{\sqrt{3}(z+\zb)(|z|^2-1)}{(1+|z|^2)^2} \theta^3 - \tfrac{\sqrt{3}(z^2+\zb^2)}{(1+|z|^2)^2} \theta^4 - \tfrac{i\sqrt{3}(z-\zb)(|z|^2-1)}{(1+|z|^2)^2} \theta^5.
\]

Since there exist two \emph{commuting} complex structures $J_\pm$ on
every 4-dimensional horizontal subspace $\calH_x^u$, the
complexification of this subspace decomposes onto the common
eigenspaces of $J_\pm$. Explicitly we have 
\[
(\calH_x^u)^\bbC = N_1\oplus N_2 \oplus \overline{N_1} \oplus \overline{N_2},
\]
where the spaces $N_1$ and $N_2$ are defined by 
$$
J_\pm N_1= i N_1,\quad\quad J_\pm N_2=\pm i N_2,$$
and $\overline{N_1}$, $\overline{N_2}$ denote their respective 
complex conjugates. The explicit formulae for the
$\tilde{g}$-duals $\tilde{n}_1$ and $\tilde{n}_2$ of the 
vectors $n_1$ and $n_2$ generating the subspaces $N_1$ and $N_2$ are the following
\begin{align*}
\tilde{n}_1 &= \tfrac{i2\sqrt{3} z (|z|^2-1)}{(1+|z|^2)^2} \theta^1 - \tfrac{2 (z^3+\zb)}{(1+|z|^2)^2} \theta^2 - \tfrac{i (1-3 z^2 -3 z\zb +z^3\zb)}{(1+|z|^2)^2} \theta^3 - \tfrac{2i (z^3-\zb)}{(1+|z|^2)^2} \theta^4 - \tfrac{1+3 z^2 -3 z\zb -z^3\zb}{(1+|z|^2)^2} \theta^5,\\
\tilde{n}_2 &= \tfrac{i2\sqrt{3} z^2}{(1+|z|^2)^2} \theta^1 + \tfrac{z^4-1}{(1+|z|^2)^2} \theta^2 - \tfrac{2iz (z^2-1)}{(1+|z|^2)^2} \theta^3 + \tfrac{i(z^4+1)}{(1+|z|^2)^2} \theta^4 + \tfrac{2z (z^2+1)}{(1+|z|^2)^2} \theta^5.
\end{align*}

The space $\mathcal{N}_\twist^{(1,0)}$ of $(1,0)$-forms with respect to the almost complex structure $\mathcal{J}_0$ is spanned by
\[
\mathcal{N}_\twist^{(1,0)} = \mathrm{Span}_\bbC (\tilde{h},\tilde{n}_1,\tilde{n}_2).
\]
Thus the integrability conditions for the CR structure $(\mathcal{N}_\twist,\mathcal{J}_0)$ have the form
\begin{align*}
\dr \tilde{u}\dz \tilde{u}\dz \tilde{h}
\dz \tilde{n}_1\dz \tilde{n}_2 &\equiv 0\\
\dr \tilde{h}\dz \tilde{u}\dz \tilde{h}\dz \tilde{n}_1\dz \tilde{n}_2 &\equiv 0\\
\dr \tilde{n}_1\dz \tilde{u}\dz \tilde{h}\dz \tilde{n}_1\dz \tilde{n}_2 &\equiv 0\\
\dr \tilde{n}_2\dz \tilde{u}\dz \tilde{h}\dz \tilde{n}_1\dz \tilde{n}_2 &\equiv 0,
\end{align*}
The expression for the other almost CR structures are analogous. 

The remaining part of proof of the theorem is skipped due to its purely computational character.\\
\end{proof}

\newcommand{\vth}{\vartheta}
\begin{remark}
We close this section with a remark that on $\twist$ there exist also
other natural geometries whose integrability conditions may encode the
torsion/curvature properties of $\sog(3)$ structures. Let us define the following real 1-forms
\begin{gather*}
  \vth^1 = \mathrm{Re} (\tilde{n}_1),\quad  \vth^2 = \mathrm{Im} (\tilde{n}_1),\quad \vth^3 = \mathrm{Re} (\tilde{n}_2),\quad \vth^4 = \mathrm{Im} (\tilde{n}_1),\\
  \vth^5 = \tilde{u},\quad \vth^6 = -\mathrm{Im} (\tilde{h}),\quad \vth^7 = \mathrm{Re} (\tilde{h}).
\end{gather*}
They define the $\tilde{g}$-orthonormal (local) coframe on $\twist$. The
 following 3-forms
\begin{align*}
\phi_1&=\tfrac{i}{2}(\tilde{n}_1\dz\overline{\tilde{n}}_1-\tilde{n}_2\dz\overline{\tilde{n}}_2)\dz\tilde{u}\\
\phi_2&=\tfrac{i}{2}(\tilde{n}_1\dz\overline{\tilde{n}}_2\dz\tilde{h}-\overline{\tilde{n}}_1\dz\tilde{n}_2\dz\overline{\tilde{h}})\\
\phi_3&=\tfrac{i}{2}\tilde{u}\dz\tilde{h}\dz\overline{\tilde{h}}
\end{align*}
are well defined on $\twist$. They may be collected to a single well
defined 3-form 
$$\phi=\phi_1+\phi_2+\phi_3.$$ This, when expressed in terms of the
orthonormal coframe $(\vth^1,\vth^2,\vth^3,\vth^4,\vth^5,\vth^6,\vth^7)$, reads 
\[
\phi=(\vth^1\wedge\vth^2 - \vth^3\wedge\vth^4)\wedge\vth^5 + (\vth^1\wedge\vth^3 - \vth^4\wedge\vth^2)\wedge\vth^6 + (\vth^1\wedge\vth^4 - \vth^2\wedge\vth^3)\wedge\vth^7 + 
\vth^5\wedge\vth^6\wedge\vth^7.
\]
It equips $\twist$ with a $G_2\subset\sog(\tilde{g})$ structure (see \cite{sal}).
\end{remark}

\section{Acknowledgements}
This paper is inspired by a talk ``Fast-hermitesche Mannigfaltigkeiten
mit paralleler charakteristischer Torsion'' which Thomas Friedrich
gave at Humboldt University on 18.05.2004. We thank Ilka Agricola, 
Thomas Friedrich, Paul-Andi Nagy and Simon Salamon for helpful
discusions. Our special thanks go to Ilka Agricola. Without her our 
collaboration would not be possible.

\end{document}

%% file: hhex.pstex_t
\begin{picture}(0,0)%
\includegraphics{hhex.pstex}%
\end{picture}%
\setlength{\unitlength}{4144sp}%
\begingroup\makeatletter\ifx\SetFigFont\undefined%
\gdef\SetFigFont#1#2#3#4#5{%
  \reset@font\fontsize{#1}{#2pt}%
  \fontfamily{#3}\fontseries{#4}\fontshape{#5}%
  \selectfont}%
\fi\endgroup%
\begin{picture}(5424,6324)(889,-5023)
\put(3826,-241){\rotatebox{50.0}{\makebox(0,0)[lb]{\smash{{\SetFigFont{12}{14.4}{\rmdefault}{\mddefault}{\updefault}{\color[rgb]{0,0,0}$\sog(3)\times\sog(1,2)$}%
}}}}}
\put(6031,-1996){\makebox(0,0)[lb]{\smash{{\SetFigFont{12}{14.4}{\rmdefault}{\mddefault}{\updefault}{\color[rgb]{0,0,0}$h_2$}%
}}}}
\put(3601,1064){\makebox(0,0)[lb]{\smash{{\SetFigFont{12}{14.4}{\rmdefault}{\mddefault}{\updefault}{\color[rgb]{0,0,0}$h_1$}%
}}}}
\put(4231,-601){\rotatebox{65.0}{\makebox(0,0)[lb]{\smash{{\SetFigFont{12}{14.4}{\rmdefault}{\mddefault}{\updefault}{\color[rgb]{0,0,0}$h_1=2 h_2$}%
}}}}}
\put(4366,-871){\makebox(0,0)[lb]{\smash{{\SetFigFont{12}{14.4}{\rmdefault}{\mddefault}{\updefault}{\color[rgb]{0,0,0}$\boxed{G_\sigma}$}%
}}}}
\put(3961,-2266){\makebox(0,0)[lb]{\smash{{\SetFigFont{12}{14.4}{\rmdefault}{\mddefault}{\updefault}{\color[rgb]{0,0,0}$\boxed{(\sog(2)\rtimes\bbR^2)\times\sog(3)}$}%
}}}}
\put(1306,-3526){\rotatebox{28.0}{\makebox(0,0)[lb]{\smash{{\SetFigFont{12}{14.4}{\familydefault}{\mddefault}{\updefault}{\color[rgb]{0,0,0}$\boxed{T\in\bgw^2_3}$}%
}}}}}
\put(1756,839){\rotatebox{297.0}{\makebox(0,0)[lb]{\smash{{\SetFigFont{12}{14.4}{\familydefault}{\mddefault}{\updefault}{\color[rgb]{0,0,0}$\boxed{T\in\bgw^2_7}$}%
}}}}}
\put(5041,-196){\rotatebox{305.0}{\makebox(0,0)[lb]{\smash{{\SetFigFont{12}{14.4}{\rmdefault}{\mddefault}{\updefault}{\color[rgb]{0,0,0}$\sog(1,2)\times\sog(3)$}%
}}}}}
\put(2408,-4732){\rotatebox{45.0}{\makebox(0,0)[lb]{\smash{{\SetFigFont{12}{14.4}{\rmdefault}{\mddefault}{\updefault}{\color[rgb]{0,0,0}$\sog(3)\times\sog(1,2)$}%
}}}}}
\put(1486,-2131){\rotatebox{300.0}{\makebox(0,0)[lb]{\smash{{\SetFigFont{12}{14.4}{\rmdefault}{\mddefault}{\updefault}{\color[rgb]{0,0,0}$\sog(1,2)\times\sog(3)$}%
}}}}}
\put(1711,-1681){\makebox(0,0)[lb]{\smash{{\SetFigFont{12}{14.4}{\familydefault}{\mddefault}{\updefault}{\color[rgb]{0,0,0}$\boxed{\bbR\times(\sog(2)\rtimes\bbR^4)}$}%
}}}}
\put(4141,-1546){\rotatebox{28.0}{\makebox(0,0)[lb]{\smash{{\SetFigFont{12}{14.4}{\rmdefault}{\mddefault}{\updefault}{\color[rgb]{0,0,0}$h_2=2 h_1$}%
}}}}}
\put(4411,-4561){\rotatebox{30.0}{\makebox(0,0)[lb]{\smash{{\SetFigFont{12}{14.4}{\rmdefault}{\mddefault}{\updefault}{\color[rgb]{0,0,0}$\sog(3)\times\sog(3)$}%
}}}}}
\put(2251,839){\rotatebox{297.0}{\makebox(0,0)[lb]{\smash{{\SetFigFont{12}{14.4}{\familydefault}{\mddefault}{\updefault}{\color[rgb]{0,0,0}$h_1=\minus 2 h_2$}%
}}}}}
\put(1981,-691){\rotatebox{25.0}{\makebox(0,0)[lb]{\smash{{\SetFigFont{12}{14.4}{\rmdefault}{\mddefault}{\updefault}{\color[rgb]{0,0,0}$\sog(3)\times\sog(3)$}%
}}}}}
\put(3871,-2896){\rotatebox{270.0}{\makebox(0,0)[lb]{\smash{{\SetFigFont{12}{14.4}{\rmdefault}{\mddefault}{\updefault}{\color[rgb]{0,0,0}$\boxed{(\sog(2)\rtimes\bbR^2)\times\sog(3)}$}%
}}}}}
\end{picture}%

%% file: hh2ex.pstex_t
\begin{picture}(0,0)%
\includegraphics{hh2ex.pstex}%
\end{picture}%
\setlength{\unitlength}{4144sp}%
\begingroup\makeatletter\ifx\SetFigFont\undefined%
\gdef\SetFigFont#1#2#3#4#5{%
  \reset@font\fontsize{#1}{#2pt}%
  \fontfamily{#3}\fontseries{#4}\fontshape{#5}%
  \selectfont}%
\fi\endgroup%
\begin{picture}(5424,6324)(889,-5023)
\put(6031,-1996){\makebox(0,0)[lb]{\smash{{\SetFigFont{12}{14.4}{\rmdefault}{\mddefault}{\updefault}{\color[rgb]{0,0,0}$h_2$}%
}}}}
\put(4726,-1276){\rotatebox{28.0}{\makebox(0,0)[lb]{\smash{{\SetFigFont{12}{14.4}{\rmdefault}{\mddefault}{\updefault}{\color[rgb]{0,0,0}$h_2=2 h_1$}%
}}}}}
\put(3601,1064){\makebox(0,0)[lb]{\smash{{\SetFigFont{12}{14.4}{\rmdefault}{\mddefault}{\updefault}{\color[rgb]{0,0,0}$h_1$}%
}}}}
\put(2431,479){\rotatebox{297.0}{\makebox(0,0)[lb]{\smash{{\SetFigFont{12}{14.4}{\familydefault}{\mddefault}{\updefault}{\color[rgb]{0,0,0}$h_1=\minus 2 h_2$}%
}}}}}
\put(4906,-871){\rotatebox{37.0}{\makebox(0,0)[lb]{\smash{{\SetFigFont{12}{14.4}{\rmdefault}{\mddefault}{\updefault}{\color[rgb]{0,0,0}$h_2=3/2 h_1$}%
}}}}}
\put(4321,-466){\rotatebox{65.0}{\makebox(0,0)[lb]{\smash{{\SetFigFont{12}{14.4}{\rmdefault}{\mddefault}{\updefault}{\color[rgb]{0,0,0}$h_1=2 h_2$}%
}}}}}
\put(1711,884){\rotatebox{300.0}{\makebox(0,0)[lb]{\smash{{\SetFigFont{12}{14.4}{\familydefault}{\mddefault}{\updefault}{\color[rgb]{0,0,0}$\boxed{T\in\bgw^2_7}$}%
}}}}}
\put(4996,-16){\rotatebox{53.0}{\makebox(0,0)[lb]{\smash{{\SetFigFont{12}{14.4}{\rmdefault}{\mddefault}{\updefault}{\color[rgb]{0,0,0}$h_1=h_2$}%
}}}}}
\put(1711,-3796){\rotatebox{45.0}{\makebox(0,0)[lb]{\smash{{\SetFigFont{12}{14.4}{\familydefault}{\mddefault}{\updefault}{\color[rgb]{0,0,0}$\boxed{K_{\odot^2_9}=0}$}%
}}}}}
\put(5311,-61){\rotatebox{53.0}{\makebox(0,0)[lb]{\smash{{\SetFigFont{12}{14.4}{\rmdefault}{\mddefault}{\updefault}{\color[rgb]{0,0,0}$\boxed{G_0}$}%
}}}}}
\put(1576,-2221){\makebox(0,0)[lb]{\smash{{\SetFigFont{12}{14.4}{\familydefault}{\mddefault}{\updefault}{\color[rgb]{0,0,0}$\boxed{K_{\wedge^1_5}=0}$}%
}}}}
\put(5356,-1411){\rotatebox{28.0}{\makebox(0,0)[lb]{\smash{{\SetFigFont{12}{14.4}{\familydefault}{\mddefault}{\updefault}{\color[rgb]{0,0,0}$\boxed{T\in\bgw^2_3}$}%
}}}}}
\put(4861,-3031){\makebox(0,0)[lb]{\smash{{\SetFigFont{12}{14.4}{\rmdefault}{\mddefault}{\updefault}{\color[rgb]{0,0,0}$\boldsymbol{G_\sigma}$}%
}}}}
\end{picture}%